\documentclass[notitlepage,11pt,reqno]{amsart}
\usepackage{amssymb,amsmath,amscd,xspace,amsthm,tocvsec2,mathrsfs}
\usepackage[breaklinks=true]{hyperref}

\usepackage[letterpaper]{geometry}
\newcommand{\To}{\ensuremath{\Rightarrow}\xspace}
\newcommand{\aaa}{\ensuremath{\alpha}\xspace}
\newcommand{\g}{\ensuremath{\gamma}\xspace}

\newcommand{\bfg}{\ensuremath{{\bf g}}\xspace}
\newcommand{\bfh}{\ensuremath{{\bf h}}\xspace}
\newcommand{\bbr}{\ensuremath{\beta_{\bf r}}\xspace}
\newcommand{\dd}{\ensuremath{\Delta}\xspace}
\newcommand{\la}{\ensuremath{\lambda}\xspace}
\newcommand{\La}{\ensuremath{\Lambda}\xspace}
\newcommand{\lal}{\ensuremath{\Lambda^\lambda_L}\xspace}
\newcommand{\sg}{\ensuremath{\Sigma_\Pi}\xspace}
\newcommand{\sgp}{\ensuremath{\Sigma_{\Pi_p}}\xspace}
\newcommand{\sgphi}{\ensuremath{\Sigma_{\Phi'}}\xspace}
\newcommand{\sgphip}{\ensuremath{\Sigma_{\Phi_p}}\xspace}

\newcommand{\sla}{\ensuremath{\Sigma_{\langle \lambda \rangle}}\xspace}
\newcommand{\ggh}{\ensuremath{G_\Pi(H)}\xspace}
\newcommand{\noti}{\ensuremath{[n] \setminus I}\xspace}
\newcommand{\notj}{\ensuremath{[n] \setminus J}\xspace}
\newcommand{\vi}{\ensuremath{\varepsilon}\xspace}
\newcommand{\B}{\ensuremath{\mathcal{B}}\xspace}
\newcommand{\calc}{\ensuremath{\mathcal{C}}\xspace}
\newcommand{\calt}{\ensuremath{\mathcal{T}}\xspace}

\newcommand{\scrh}{\ensuremath{\mathscr{H}}\xspace}

\newcommand{\N}{\ensuremath{\mathbb{N}}\xspace}
\newcommand{\Z}{\ensuremath{\mathbb{Z}}\xspace}
\newcommand{\Q}{\ensuremath{\mathbb{Q}}\xspace}
\newcommand{\comp}{\ensuremath{\mathbb{C}}\xspace}

\newcommand{\sql}{\ensuremath{\sqrt[l]{1}}\xspace}
\newcommand{\uql}{\ensuremath{u_{q,l}}\xspace}
\newcommand{\mfh}{\ensuremath{\mathfrak{h}}\xspace}
\newcommand{\mfg}{\ensuremath{\mathfrak{g}}\xspace}
\newcommand{\mfgone}{\ensuremath{\mathfrak{g}_{ab}}\xspace}
\newcommand{\mfb}{\ensuremath{\mathfrak{b}}\xspace}
\newcommand{\mfgs}{\ensuremath{\mathfrak{G}^*}\xspace}
\newcommand{\mfgn}[1]{\ensuremath{\mathfrak{G}_{n,#1}}\xspace}
\newcommand{\aqn}[1]{\ensuremath{\mathcal{A}_{q,#1}(n)}\xspace}
\newcommand{\one}[1]{\ensuremath{{#1}_{(1)}}\xspace}
\newcommand{\two}[1]{\ensuremath{{#1}_{(2)}}\xspace}
\newcommand{\three}[1]{\ensuremath{{#1}_{(3)}}\xspace}
\newcommand{\tangle}[1]{\ensuremath{\langle #1 \rangle}\xspace}
\newcommand{\nomial}[1]{\ensuremath{\binom{#1}{l_1, \dots, l_{k-1} }}\xspace}
\newcommand{\comment}[1]{}

\DeclareMathOperator{\id}{\ensuremath{id}\xspace}
\DeclareMathOperator{\ad}{\ensuremath{ad}\xspace}
\DeclareMathOperator{\ch}{\ensuremath{char}\xspace}
\DeclareMathOperator{\im}{\ensuremath{im}\xspace}
\DeclareMathOperator{\Hom}{\ensuremath{Hom}\xspace}
\DeclareMathOperator{\Sym}{\ensuremath{Sym}\xspace}
\DeclareMathOperator{\perm}{\ensuremath{perm}\xspace}

\theoremstyle{plain}
\newtheorem{theorem}[equation]{Theorem}
\newtheorem{lemma}[equation]{Lemma}
\newtheorem{prop}[equation]{Proposition}
\newtheorem{cor}[equation]{Corollary}

\theoremstyle{definition}
\newtheorem{defn}[equation]{Definition}
\newtheorem{stand}[equation]{Standing Assumption}
\newtheorem{remark}[equation]{Remark}
\newtheorem{example}[equation]{Example}

\numberwithin{equation}{section}

\begin{document}
\vspace*{-0.5ex}
\title{The sum of a finite group of weights of a Hopf algebra}
\author{Apoorva Khare}
\address{Department of Mathematics and Department of Statistics, Stanford
University, 390 Serra Mall, Stanford, CA - 94305, USA}
\email{khare@stanford.edu}
\date{\today}

\begin{abstract}
Motivated by the orthogonality relations for irreducible characters of a
finite group, we evaluate the sum of a finite group of linear characters
of a Hopf algebra, at all grouplike and skew-primitive elements. We then
discuss results for products of skew-primitive elements. Examples include
groups, (quantum groups over)
Lie algebras, the small quantum groups of Lusztig, and their variations
(by Andruskiewitsch and Schneider).
\end{abstract}
\maketitle

\settocdepth{section}
{\tt \tableofcontents}

\section{Introduction}

\subsection{Motivation}

Suppose $G$ is a finite group, with irreducible characters $\widehat{G} =
\{ \chi : \comp G \to \comp \}$, over $\comp$. As is well known, they
satisfy the {\it orthogonality relations}. Here is a consequence: if one
defines $\theta_\chi := \frac{1}{|G|} \sum_{g \in G} \chi(g^{-1}) g \in H
= \comp G$, then $\nu(\theta_\chi) = 0$ for $\chi \neq \nu \in
\widehat{G}$. Similarly, the other orthogonality relation (for columns)
implies that $\sum_{\chi \in \widehat{G}} (\dim \rho_\chi)\chi(g)$ is
either zero or a factor of $|G|$, depending on whether or not $g = 1$.
(Here, $\rho_\chi$ is the irreducible representation with character
$\chi$.)\smallskip

In what follows, we work over a commutative integral domain $R$. Note
that there is an analogue of the first orthogonality relation for any
$R$-algebra $H$ that is a free $R$-module. Namely, given a linear
character (i.e., algebra map or {\it weight}) $\la : H \to R$, define a
{\it left $\la$-integral} of $H$ to be any (nonzero) $\lal \in H$ so that
$h \lal = \la(h) \lal$ for all $h \in H$. One can similarly define right
and two-sided $\la$-integrals in $H$. (For instance, the $\theta_\chi$'s
above, are two-sided $\chi$-integrals for any weight $\chi$.) Then the
following result holds in any algebra:

\begin{lemma}
If $\la \neq \nu$ are weights of $H$ with corresponding nonzero left
integrals $\lal, \La^\nu_L$ respectively, then $\nu(\lal) = 0 = \lal
\La^\nu_L$.
\end{lemma}

\begin{proof}
Choose $h$ so that $\la(h) \neq \nu(h)$. Then
\[ \nu(h) \nu(\lal)\La^\nu_L = \nu(h) \lal \La^\nu_L = h \lal \La^\nu_L =
(h \lal) \La^\nu_L = \la(h) \lal \La^\nu_L = \la(h) \nu(\lal) \La^\nu_L,
\]

\noindent so $\nu(\lal) (\la(h) - \nu(h))\La^\nu_L = 0$. Since we are
working over an integral domain and within a free module, this implies
that $\nu(\lal) = 0$. Moreover, $\lal \La^\nu_L = \nu(\lal) \La^\nu_L =
0$.
\end{proof}

It is now natural to seek a ``Hopf-theoretic" analogue for the second
orthogonality relation (which might involve only weights, and not all
irreducible characters). Note that in general, a Hopf algebra might not
have a nontrivial (sub)group of linear characters that is finite: for
instance, $U(\mfg)$ for a complex Lie algebra $\mfg$.
However, if such subgroups do exist, then we attempt to evaluate the sum
of all weights in the subgroup (i.e., the coefficient of each weight is
$\dim \rho = 1$), at various elements of $H$. For example, one may ask:
does this sum always vanish at a nontrivial grouplike $g \in H$?\medskip

This problem is also interesting from another perspective. Given a Hopf
algebra $H$, it is interesting to seek connections between the
representations of $H$ and those of $H^*$ (or, the comodules over $H^*$
and over $H$ respectively). A famous example where both of these come
into play is the well-known formula due to Radford \cite{Rad} concerning
the antipode.
The setting of this paper is another such case: note that the weights of
$H$ are precisely the one-dimensional representations of $H^*$ (or the
group-like elements of $H^*$), while grouplike and skew-primitive
elements in a Hopf algebra $H$ correspond to the $1$-dimensional and
non-semisimple $2$-dimensional $H$-comodules (or $H^*$-modules)
respectively. The present paper explores the interplay between these
objects.

\subsection{One of the setups, and some references}

Instead of attempting a summary of the results, which are several and
computational, we make some remarks. Given an $R$-algebra $A$, let
$\Gamma_A$ denote the set of weights of $A$. (If $H$ is a Hopf algebra,
then $\Gamma_H$ is a group.) Let $\Pi \subset \Gamma_H$ be a finite group
of weights of a Hopf algebra $H$ over a commutative unital integral
domain $R$. One has the notion of {\it grouplike} elements (i.e., $\dd(g)
= g \otimes g$) and {\it skew-primitive} elements in $H$.

The first step is to compute $\sg$ at all grouplike and skew-primitive
elements in $H$, where $\sg := \sum_{\g \in \Pi} \g : H \to R$. In a wide
variety of examples - including finite-dimensional pointed Hopf algebras
\cite{AS} - the computations reduce to grouplike elements. In other
words, there are several families of algebras, where knowing $\sg$ at
grouplike elements effectively tells us $\sg$ at {\it all} elements.

The next objective is to evaluate $\sg$ at products of skew-primitive
elements. Once again, in the spirit of the previous paragraph, there are
numerous examples of Hopf algebras generated by grouplike and
skew-primitive elements in the literature. The first two examples below
are from folklore, and references can be found in \cite{Mon}.
\begin{enumerate}
\item By the Cartier-Kostant-Milnor-Moore Theorem (e.g., see
\cite[Theorem 5.6.5]{Mon}), every cocommutative connected Hopf algebra
$H$ over a field of characteristic zero, is of the form $U(\mfg)$, where
$\mfg$ is the set of primitive elements in $H$. Similarly, every complex
cocommutative Hopf algebra is generated by primitive and grouplike
elements.

\item If the Hopf algebra is pointed (and over a field), then by the
Taft-Wilson Theorem \cite[Theorem 5.4.1.1]{Mon}, our results can evaluate
$\sg$ on any element of $C_1$, the first term in the coradical filtration
(which is spanned by grouplike and skew-primitive elements).

\item The final example is from a recent paper \cite{AS}. The
Classification Theorem 0.1 says, in particular, that if $H$ is a
finite-dimensional pointed Hopf algebra over an algebraically closed
field of characteristic zero, and the grouplike elements form an abelian
group of order coprime to 210, then $H$ is generated by grouplike and
skew-primitive elements, and is a variation of a small quantum group of
Lusztig.
\end{enumerate}

We conclude this section with one of our results. Say that an element $h
\in H$ is {\it pseudo-primitive} with respect to $\Pi$ if $\dd(h) = g
\otimes h + h \otimes g'$ for grouplike $g,g'$ satisfying $\g(g) =
\g(g')$ for all $\g \in \Pi$.

\begin{theorem}
Fix $n \in\N$, as well as the $R$-Hopf algebra $H$ and a finite subgroup
of weights $\Pi$ of $H$. Suppose $h_1, \dots, h_n \in H$ are
pseudo-primitive with respect to $\Pi$, and $\dd(h_i) = g_i \otimes h_i +
h_i \otimes g'_i$ for all $i$. Define $\bfh = \prod_i h_i$, and
similarly, $\bfg, \bfg'$.
\begin{enumerate}
\item If $\ch(R) = 0$ or $\ch(R) \nmid |\Pi|$, then $\sg(\bfh) = 0$.

\item Suppose $0<p=\ch(R)$ divides $|\Pi|$, and $\Pi_p$ is any $p$-Sylow
subgroup. If $\Pi_p \ncong (\Z / p\Z)^m$ for any $m>0$, then $\sg(\bfh) =
0$.

\item ($p$ as above.) Define $\Phi := \Pi / [\Pi,\Pi]$, and by above,
suppose $\Phi_p \cong (\Z / p\Z)^k$ is a $p$-Sylow subgroup of $\Phi$.
Let $\Phi'$ be any Hall complement(ary subgroup); thus $|\Phi'| = |\Phi|
/ |\Phi_p|$. Then
\[ \sg(\bfh) = |[\Pi, \Pi]| \cdot \sgphi(\bfg) \cdot \sgphip(\bfh). \]

\item If $\sgphip(\bfh)$ is nonzero, then $(p-1)|n$, and $0 \leq k \leq
n/(p-1)$. (Moreover, examples exist wherein $\sgphip(\bfh)$ can take any
value $r \in R$.)
\end{enumerate}
\end{theorem}

\noindent These results occur below as Proposition \ref{P7}, and Theorems
\ref{T6}, \ref{T5}, and \ref{T8} respectively.

\subsection{Organization}

We quickly explain the organization of this paper. In Section \ref{S2},
we compute $\sg$ at all grouplike elements in a Hopf algebra. This turns
out to be extremely useful in computing $\sg$ in a large class of
examples. More precisely, many algebras $A$ that are well-studied in the
literature, contain a Hopf (sub)algebra $H$, and whose sets of weights
are subsets of $\Gamma_H$ and kill a large subspace of $A$.
These examples include quantum groups, the quantum Virasoro algebra, and
finite-dimensional pointed Hopf algebras.

In Section \ref{S3}, we evaluate $\sg$ at all skew-primitive elements of
$H$. This is followed by a brief remark concerning the ``degenerate"
example of quiver (co)algebras.

Computing $\sg$ at all products of skew-primitive elements is a difficult
problem. We show in the next section \ref{Ssubquot} that on occasion, it
can be reduced to computing $\Sigma_{\Gamma_H / [\Gamma_H, \Gamma_H]}$.
In Section \ref{S6}, we are able to obtain results when $\Gamma_H$ itself
is abelian. These computations are useful in working with products of
special kinds of skew-primitive elements in other sections.

In the rest of the paper, we work with ``pseudo-primitive elements" $h_i$
in order to obtain more detailed results. In Section \ref{S5}, we show
that $\sg(h_1 \dots h_n) = 0$ whenever $\ch(R) \nmid |\Pi|$. In Section
\ref{Spseudo}, we study the case when $\ch(R) = p$ divides $|\Pi|$. In
this case, there are severe restrictions on the $p$-Sylow subgroup of
$\Pi$, in order for $\sg({\bf h})$ to be nonzero; moreover, we write down
a result that helps compute $\sg({\bf h})$.

We conclude with a detailed study of further examples - Lie algebras,
degenerate affine Hecke algebras of reductive type, and then a Hopf
algebra generated by grouplike and skew-primitive elements, where
$\sg({\bf h})$ can take on all values.

\section{Grouplike elements and quantum groups}\label{S2}

\subsection{Preliminaries}

We first set some notation, and make some definitions.

\begin{defn}\label{D1}
Suppose $R$ is a commutative unital integral domain.
\begin{enumerate}
\item {\it Integers} in $R$ are the image of the group homomorphism $\Z
\to R$, sending $1 \mapsto 1$.

\item A {\it weight} of an $R$-algebra $H$ is an $R$-algebra map $: H \to
R$. Denote the set of weights by $\Gamma_H$. Occasionally we will also
use $\Gamma = \Gamma_H$. Given $\nu \in \Gamma_H$, the $\nu$-{\it weight
space} of an $H$-module $V$ is
$V_\nu := \{ v \in V : h \cdot v = \nu(h) v\ \forall h \in H \}$.

\item Given a left $R$-module $H$, define $H^* := \Hom_{R-{\rm
mod}}(H,R)$.

\item An $R$-Hopf algebra $H$ is an $R$-algebra $(H, \mu = \cdot, \eta)$
(where $\mu,\eta$ are coalgebra maps) that is also an $R$-coalgebra $(H,
\dd, \vi)$ (where $\dd, \vi$ are algebra maps), further equipped with an
antipode $S$ (which is an $R$-(co)algebra anti-homomorphism).

\item In a Hopf algebra (or a bialgebra), an element $h$ is {\it
grouplike} if $\dd(h) = h \otimes h$, and {\it primitive} if $\dd(h) = 1
\otimes h + h \otimes 1$. Define $G(H)$ (respectively $H_{prim}$) to be
the set of grouplike (respectively primitive) elements in a Hopf algebra
$H$.
\end{enumerate}
\end{defn}

There are several standard texts on Hopf algebras; for instance, see
\cite{Abe,DNR,Mon,Swe}. In particular, since $H$ is a coalgebra, $H^*$ is
also an $R$-algebra under {\it convolution} $\dd^*$: given $\la,\nu \in
H^*$ and $h \in H$, one defines $\tangle{\la * \nu, h} := \tangle{\la
\otimes \nu, \dd(h)}$. By \cite[Theorem 2.1.5]{Abe} (also see \cite[Lemma
4.0.3]{Swe}), the set $\Gamma_H$ of weights is now a group under $*$,
with inverse given by $\tangle{\la^{-1}, h} := \tangle{\la, S(h)}$.

Note that for an algebra $H$ over a field $k$, the dual space $H^*$ is
not a coalgebra in general. However, define $H^\circ$ to be the set of
linear functionals $f : H \to k$ whose kernel contains an ideal of finite
$k$-codimension. Then by \cite[Proposition 1.5.3 and Remark 1.5.9]{DNR},
$H^\circ$ is a coalgebra whose set of grouplike elements is precisely
$\Gamma_H$.

\begin{stand}
For this article, $H$ is any Hopf algebra
over a commutative unital integral domain $R$.
Fix a finite subgroup of weights $\Pi \subset \Gamma = \Gamma_H$.\medskip
\end{stand}

In general, given a finite subgroup $\Pi \subset \Gamma_H$ for any Hopf
algebra $H$, the element $\sg := \sum_{\g \in \Pi} \g$ is a functional in
$H^*$, and if $R$ (i.e., its quotient field) has characteristic zero,
then $\sg$ does not kill the scalars $\eta(R)$. What, then, is its
kernel? How about if $\Pi$ is cyclic, or all of $\Gamma_H$ (this, only if
$H$ is $R$-free, and finite-dimensional over the quotient field of $R$)?

\begin{lemma}\label{L2}
Suppose $H$ is a Hopf algebra, and $\Pi \subset \Gamma_H$ is a finite
subgroup of weights.
\begin{enumerate}
\item $\sg(1) = 0$ if and only if $\ch(R)$ divides $|\Pi|$.

\item $[H,H] \subset \ker \sg$.

\item $\sg(\ad h(h')) = \vi(h) \sg(h')$ for all $h,h' \in H$. In
particular, if $h \in \ker \vi$, then $\im( \ad h ) \subset \ker \sg$.
\end{enumerate}
\end{lemma}

\noindent Here, $\ad$ stands for the usual {\it adjoint action} of $H$.
In other words, $\ad h(h') := \sum \one{h} h' S(\two{h})$, where we use
{\it Sweedler notation}: $\Delta(h) = \sum \one{h} \otimes \two{h}$.

\begin{proof}
The first part is easy, and the other two follow because the statements
hold if $\sg$ is replaced by any (algebra map) $\g \in \Gamma_H$.
\end{proof}

The goal of this section and the next, is to evaluate $\sg$ at all
grouplike and skew-primitive elements in $H$. For these computations, a
key fact to note is that for all $\la \in \Pi$, the following holds in
the $R$-algebra $H^*$:
\begin{equation}\label{E0}
\la * \sg = \sum_{\nu \in \Pi} \la * \nu = \sg = \dots = \sg * \la.
\end{equation}



\begin{remark}
We occasionally compute $\sg$ with $H$ an $R$-algebra (that is not a Hopf
algebra), where $\Pi \subset \Gamma_H$ has a group structure on it. As
seen in Proposition \ref{Pq} below, there is an underlying Hopf algebra
in some cases.
\end{remark}

\begin{defn}
Suppose we have a subset $\Theta \subset \Gamma = \Gamma_H$, and $\la \in
\Pi$.
\begin{enumerate}
\item For $g \in G(H)$, set $\Gamma_g := \{ \g \in \Gamma : \g(g) = 1
\}$, and $\Theta_g := \Gamma_g \cap \Theta$.

\item $G_{\Theta}(H)$ is the set (actually, normal subgroup) of grouplike
elements $g \in G(H)$ so that $\g(g) = 1$ for all $\g \in \Theta$.

\item For finite $\Theta$, the functional $\Sigma_{\Theta} \in H^*$ is
given by $\Sigma_{\Theta} := \sum_{\g \in \Theta} \g$. Also set
$\Sigma_{\emptyset} := 0$.

\item $n_\la := o_\Pi(\la) = |\tangle{\la}|$ is the order of $\la$ in
$\Pi$.
\end{enumerate}
\end{defn}

\begin{remark}\hfill
\begin{enumerate}
\item For instance, $G_{\{ \vi \}}(H) = G(H)$, and $G(H) \cap [1 + \im
(\id - S^2)] \subset G_\Gamma(H)$ because every $\g \in \Gamma =
\Gamma_H$ equals $(\g^{-1})^{-1}$. This follows since from earlier in
this section, $\Gamma_H$ is a group with unit $\varepsilon$ and inverse
given by $\gamma \mapsto \gamma \circ S$.

%

\item For any $g \in G(H)$ and $\Theta \subset \Gamma_H$, $\Theta_g =
\Theta$ if and only if $g \in G_{\Theta}(H)$.

\item $\Theta_g$ is a subgroup if $\Theta$ is.
\end{enumerate}
\end{remark}

\subsection{Grouplike elements}

We first determine how $\sg$ acts on grouplike elements, and answer the
motivating question above, of finding a Hopf-theoretic analogue of the
second orthogonality relations for group characters.

\begin{prop}[``Orthogonality" at grouplike elements]\label{P5}
If $g \in G(H) \setminus \ggh$, then $\sg(g) = 0$. If $g \in \ggh$, then
$\sg(g) = |\Pi|$.
\end{prop}

\begin{proof}
To show the first part, apply equation \eqref{E0} to $g$, with $\la
\notin \Gamma_g$. The second part is obvious.
\comment{
If $g \in \ggh$, the result is trivial, so now assume that $g \notin
\ggh$. Then there is some weight $\la$ that does not send it to 1. By
equation \eqref{E0}, we get that $\sg(g) = (\la * \sg)(g) = \la(g)
\sg(g)$. Thus $(1 - \la(g)) \sg(g) = 0$, and we are done since $\la(g)
\neq 1$.

Alternatively now, $\Pi$ is a finite group, and it has a subgroup
$\tangle{\la}$ generated by $\la$. Thus $\la(h)$ is a nontrivial root of
unity, so $\sla(h) = \sum_{\g \in \tangle{\la}} \g(h) = 0$. If $\B$ is a
set of coset representatives of $\tangle{\la}$ in $\Pi$, then $\Pi = \B *
\tangle{\la}$, so
\[ \sum_{\nu \in \Pi} \nu(h) =  \sum_{\beta \in \B,\ \g \in \tangle{\la}}
(\beta * \g)(h) = \sum_{\B * \tangle{\la}} (\beta \otimes \g)(\dd(h)) =
\sum_{\beta \in \B} \beta(h) \cdot \sla(h) = 0 \]

say $\{ \vi = 1_\Gamma, \la, \dots, \la^{n-1} \}$. We first note that for
any $\la, \g \in \Gamma$,
\[ (\la * \g)(h) = (\la \otimes \g)(\dd(h)) = \la(h) \g(h) \]

\noindent whence $\la^n = \vi$ implies that $\la(h)^n = \vi(h) = 1$.
Hence
\[ \sum_{i=0}^{n-1} \la^i(h) = \sum_{i=0}^{n-1} \la(h)^i = \frac{\la(h)^n
- 1}{\la(h)-1} = 0 \]

\noindent since $\la(h) \neq 1$.

But now, we sum over a set $\B$ of coset representatives of the
subgroup $\tangle{\la}$ inside $\Gamma$, to get
\[ \sum_{\g \in \Gamma} \g(h) = \sum_{\g \in \B} \sum_{i=0}^{n-1} (\g *
\la^i)(h) = \sum_{\B} \g(h) \sum_{i=0}^{n-1} \la^i(h) = 0 \]
 
\noindent and we are done.
}
\end{proof}

\noindent We now introduce some notation, which is used in discussing
several examples.

\begin{defn}
Let $R$ be a commutative unital integral domain, and $l$ be a nonnegative
integer.
\begin{enumerate}
\item $\sqrt{1}$ (respectively $\sql$) is the set of ($l^{\rm th}$) roots
of unity in $R$. (Thus, $\sqrt[0]{1} = R^\times,\ \sqrt[1]{1}= \{ 1 \},$
and $\sqrt{1} = \cup_{l>0} \sql$.)

\item Given $q \in R^\times$, $\ch(q)$ is the smallest positive integer
$m$ so that $q^m = 1$, and zero if no such $m$ exists.

\item The group $\mfgn{l}$ is the abelian group generated by $\{ K_i : 1
\leq i \leq n \}$, with relations $K_i K_j = K_j K_i, K_i^l = 1$.

\item The group $\mfgs_{n,l}$ is defined to be $(\sqrt{1} \cap \sql)^n$.
\end{enumerate}
\end{defn}

\noindent Thus, $\mfgn{l}$ is free if and only if $l=0,\ \mfgs_{n,l} =
(\sql)^n\ \forall l>0$, and $\mfgs_{n,0} = (\sqrt{1})^n$.\medskip

\noindent In light of the motivation, the first example where we apply
the above result, is:

\begin{example}[Group rings]\label{Egroup}
The above result computes $\sg$ on all of $H$, if $H$ is a group ring. We
present a specific example: $G = \mfgn{l}$ (defined above), for (a fixed)
$n \in \N$ and $l \geq 0$.
Then $\Gamma = (\sql)^n$, and any finite order element $\g \in \Gamma$
maps each $K_i$ to a root of unity in $R$. Thus, $\Pi \subset
(\sqrt{1})^n \cap \Gamma = \mfgs_{n,l}$.

We now compute $\sg(g)$ for some $g = \prod_{i=1}^n K_i^{n_i}$, where
$n_i \in \Z\ \forall i$. Note that the set $\{ \g(K_i) : 1 \leq i \leq
n,\ \g \in \Pi \}$ is a finite set of roots of unity; hence the subgroup
of $\sqrt{1} \subset R^\times$ that it generates, is cyclic, say
$\tangle{\zeta}$. Thus, $\g(K_i) = \zeta^{l_i(\g)}$ for some $l_i : \Pi
\to \Z$.

The above result now says that $\sg(g) = 0$ if there exists $\g \in \Pi$
so that $|\tangle{\zeta}|$ does not divide $\sum_{i=1}^n n_i l_i(\g)$,
and $|\Pi|$ otherwise. (Of course, one can also apply the above result
directly to $g = \prod_i K_i^{n_i}$.)
\end{example}

\begin{remark}
From the above proposition, finding out if $g \in \ggh$ is an important
step. However, since $\beta(g) = 1\ \forall \beta \in [\Pi, \Pi]$, hence
it suffices to compute if $\la(g) = 1$, where the $\la$'s are the lifts
of a set of generators of the (finite) abelian group $\Pi / [\Pi, \Pi]$.
We see more on this in Section \ref{Ssubquot} below.
\end{remark}

\subsection{Application to quantum groups and related examples}

We now mention some more examples where the above result applies: Hopf
algebras that quantize semisimple Lie algebras, their Borel subalgebras,
and polynomial algebras (i.e., coordinate rings of affine spaces/abelian
Lie algebras). There are yet other algebras mentioned below, which are
not Hopf algebras but can be treated similarly.

To discuss these examples, some more basic results are needed; here is
the setup for them. Suppose an $R$-algebra $A$ contains a Hopf subalgebra
$H$, so that $A$ is an $\ad H$-module (with possible weight spaces
$A_\nu$). Then one has the following result:

\begin{lemma}\label{Lkey}
Suppose $\mu \in \Gamma_A$ (i.e., $\mu : A \to R$ is an $R$-algebra map).
If $\nu \in \Gamma_H$ and $\nu \neq \vi$, then $\mu \equiv 0$ on $A_\nu$.
\end{lemma}

\begin{proof}
Given $\nu \neq \vi$, choose $h \in H$ such that $\nu(h) \neq \vi(h)$.
Now given $a_\nu \in A_\nu$, apply $\mu$ to the equation: $\ad h(a_\nu) =
\nu(h) a_\nu$. Simplifying this yields: $\vi(h) \mu(a_\nu) = \nu(h)
\mu(a_\nu)$, and since $R$ is an  integral domain, $\mu(a_\nu) = 0$.
\end{proof}

\noindent Applying this easily yields the following result.

\begin{prop}\label{Pq}
($A,H$ as above.) Suppose an $R$-algebra $A$ contains $H$ and a vector
subspace $V$, that is of the form $V = \oplus_{\nu \neq \vi} V_\nu$ (for
the $\ad H$-action).
\begin{enumerate}
\item Every $\mu \in \Gamma_A$ kills $AVA$.

\item If $A = H + AVA$, then $\Gamma_A \subset \Gamma_H$.

\item Say $A = H + AVA$, and $\Pi \subset \Gamma_A$ is a finite subgroup
of weights of $A$ (from above). If $a \in A$ satisfies $a - \sum_{g \in
G(H)} a_g g \in AVA$ (where $a_g \in R\ \forall g$), then $\sg(a) = |\Pi|
\sum_{g \in \ggh} a_g$.
\end{enumerate}
\end{prop}

\begin{remark}
Thus, if $\Pi \subset \Gamma_A$ is a group (i.e., a subset with a group
structure on it), then $\sg(AVA) = 0$. Hence, computing $\sg$ at any $a
\in A$ essentially reduces to the case of the Hopf subalgebra $H$. When
$H$ is a group algebra, Proposition \ref{P5} above tells us the answer in
this case - assuming that the group operation in $\Pi$ agrees with the
one in $\Gamma_H$.

Moreover, even though $A$ is not a Hopf algebra here, note that the
computations come from Hopf algebra calculations (for $H$).
\end{remark}

\begin{proof}
The first part follows from Lemma \ref{Lkey}, the third part now follows
from Proposition \ref{P5} and Lemma \ref{Lkey}, and the second part
follows by observing that $\mu : H \to R$ is an algebra map only if
$\mu|_G \in \Gamma_G$ and $\mu|_V \equiv 0$. (Additional relations in $H$
may prevent every $\mu \in \Gamma_H$ from being a weight in $\Gamma_A$.)
\end{proof}

\noindent It is now possible to apply the above theory to some examples;
note that they are not always Hopf algebras. In each case, $G$ is of the
form $\mfgn{l}$ for some $n,l$. Also choose a special element $q \in
R^\times$ in each case; then $\ch(q) | l$.

\begin{example}[``Restricted" quantum groups of semisimple Lie algebras]
For this example, $R = k$ is a field with $\ch k \neq 2$, with a special
element $q \neq 0,\pm 1$. Suppose $\mfg$ is a semisimple Lie algebra over
$k$, together with a fixed Cartan subalgebra and root space decomposition
(e.g., using a Chevalley basis, as in \cite[Chapter 7]{Hum1}).

One then defines the (Hopf) algebra $U_q(\mfg)$ as in \cite[\S
4.2,4.3]{Ja2}. In particular, note that it is generated by $\{ K_j^{\pm
1}, e_j, f_j : 1 \leq j \leq n \}$ (here, $n$ is the rank of $\mfg$, and
the $\aaa_i$ are simple roots), modulo the relations:
\[ K_i e_j K_i^{-1} = q^{(\aaa_i, \aaa_j)} e_j,\
K_i f_j K_i^{-1} = q^{-(\aaa_i, \aaa_j)} f_j,\quad
e_i f_j - f_j e_i = \delta_{ij} \frac{K_i - K_i^{-1}}{q_i - q_i^{-1}}, \]

\noindent where $q_i = q^{(\aaa_i, \aaa_i)/2}$ for some bilinear form
$(\cdot, \cdot)$ on $\mfh^*$. (We may also need that $q^4$ or $q^6$ is
not 1, and possibly also that $\ch k \neq 3$.) The other relations are
that $K_i K_j = K_j K_i,\ K_i^{\pm 1} K_i^{\mp 1} = 1$, and the (two)
quantum Serre relations. Define $V := \oplus_{j=1}^n (ke_j \oplus kf_j)$.

Now define the ``restricted" quantum group as in \cite[Chapter 6]{Maj}.
More precisely, given some fixed $l \geq 0$ so that $q^l = 1$ (whence
$\ch(q) | l$), define the associative (not necessarily Hopf) algebra
$\uql(\mfg)$ to be the quotient of $U_q(\mfg)$ by the relations (for all
$j$) $K_j^l = 1$, and $e_j^l = f_j^l = 0$ if $l>0$. Note that $\uql(\mfg)
= U_q(\mfg)$ if $l=0$, and $\uql(\mfg)$ is a Hopf algebra if $l=0$ or
$\ch(q)$. Moreover, Proposition \ref{Pq} allows us to compute $\sg$ for
all $l$.\medskip

For each $j$, note that $e_j,f_j$ are weight vectors (with respect to the
adjoint action of the abelian group $G := \mfgn{l}$ generated by all
$K_i^{\pm 1}$) with weights $q^{\pm \aaa_j} \neq \vi = q^0$. Hence Lemma
\ref{Lkey} implies that $\mu(e_j) = \mu(f_j) = 0\ \forall j, \mu$.
Moreover, given the PBW property for $A = U_q(\mfg)$, we know that $A =
k\mfgn{l} \bigoplus (V_- A + A V_+)$, where $V_+, V_-$ are the spans of
the $e_j$'s and $f_j$'s respectively. Hence $\mu \in \Gamma_A \subset
\Gamma_G$ by the result above.

Every $\mu \in \Gamma_G$ is compatible with the commuting of the $K_i$'s, the
quantum Serre relations (since $\mu|_V \equiv 0$), and the ``$l^{\rm th}$
power relations". The only restriction is the last one left, namely: $0 =
\mu([e_i, f_i])$, which gives us that $\mu(K_i) = \mu(K_i^{-1})$ for all
$i, \mu$. Hence $\mu(K_i) = \pm 1$, so that $\Gamma_A \cong (\Z / 2\Z)^n
\cap (\sql)^n$, which is of size $2^n$ or 1, depending on whether $l$ is
even or odd. We now compute $\sg(a)$ using the second part of Proposition
\ref{Pq}, for any $a \in A$.
\end{example}

\begin{example}[Restricted quantum groups of Borel subalgebras] We
consider the subalgebra $A' = \uql(\mfb)$ of $A = \uql(\mfg)$, that is
generated by $\{ K_i^{\pm 1}, e_i : 1 \leq i \leq n \}$; once again, this
algebra quantizes the Borel subalgebra $\mfb$ of $\mfg$ if $l=0$ (and is
a Hopf algebra if $l=0$ or $\ch(q)$). Moreover, every $\mu \in \Gamma_{A'}$
kills each $e_i$ (where we use $V_+$ for $V$), and we have $\Gamma_{A'} \subset
\Gamma_G$, by Proposition \ref{Pq} above. Moreover, all such maps $\mu
\in \Gamma_G$ are admissible (i.e., extend to all of $A'$), so
$\Gamma_{A'} = \Gamma_G \cong (\sql)^n$.

Now if $\Pi \subset \Gamma_{A'}$ is a finite subgroup, then as above,
$\Pi \subset \mfgs_{n,l}$ (note that $R=k$ here), and furthermore,
evaluation of $\sg$ once again reduces to the grouplike case.
\end{example}

\begin{example}[Taft algebras]
Given a primitive $n$th root $q$ of unity, the {\it $n$th Taft algebra}
is
\[ T_n := R \langle x,g \rangle / (g x - q x g,\ g^n - 1,\ x^n). \]

\noindent Once again, every weight must kill $x$, and sends $g$ to some
power of $q$. Hence the set $\Gamma_{T_n}$ of weights is cyclic, whence
so is $\Pi$. It is now easy to show:

\begin{lemma}
Every weight kills $x$, and $\sg(g^k) = |\Pi|$ if $|\Pi|$ divides $k$, or
$0$ otherwise.
\end{lemma}
\end{example}

\begin{example}[Quantization of affine space]
We refer to \cite{Hu2}; once again, $R$ is a unital commutative integral
domain. The quantum affine space over $R$ (with a fixed element $q \in
R^\times$) is the quadratic algebra $T_R(V) / (x_j x_i - q x_i x_j, i <
j)$, where $V := \oplus_i Rx_i$. This does not have a Hopf algebra
structure; however, Hu presents a quantization of $R[x] := \Sym_R(V)$ in
\cite[\S 5]{Hu2} - that is, the quantum group associated to a
``finite-dimensional" abelian Lie algebra.
Now consider a more general associative (not necessarily Hopf)
$R$-algebra $A = \aqn{l}$ generated by $\{ K_i^{\pm 1}, x_i : 1 \leq i
\leq n \}$, with the relations:
\begin{eqnarray*}
&& K_i K_j = K_j K_i,\qquad K_i^l = 1,\qquad K_i^{\pm 1} K_i^{\mp 1} =
1,\\
&& K_i x_j K_i^{-1} = \theta_{ij} q^{\delta_{ij}} x_j,\qquad x_i x_j =
\theta_{ij} x_j x_i,
\end{eqnarray*}

\noindent where $\ch(q) | l$, and $\theta_{ij}$ equals $q$ (respectively
$1, q^{-1}$) if $i > j$ (respectively $i=j,\ i<j$). Note that
$\aqn{\ch(q)} = \mathcal{A}_q(n)$, the Hopf algebra introduced and
studied by Hu, and $\aqn{l}$ becomes the Hopf algebra $R[x]$ if
$q=l=1$.\medskip

We consider the ``nontrivial" case $q \neq 1$ (the $q=1$ case is
discussed later). Once again, each $x_j$ is a weight vector with respect
to the (free abelian) group $G = \mfgn{l}$, and no $x_j$ is in the
$\vi$-weight space, so every $\mu \in \Gamma_A$ kills $x_j$ for all $j$.
As in the previous example, $\Gamma_A \cong (\sql)^n$, and any finite
subgroup $\Pi$ must be contained in $\mfgs_{n,l}$.

Moreover, evaluation of $\sg$ once again reduces to the grouplike
case.
\end{example}

\begin{example}[Quantization of the Virasoro algebra] For this example,
assume that $R$ is a field, and $q \in R^\times$ is not a root of unity.
Now refer to \cite[Page 100]{Hu1} for the definitions; the Hopf algebra
in question is the $R$-algebra $\mathcal{U}_q$ generated by $\calt,
\calt^{-1}, c, e_m (m \in \Z)$ with relations:
\begin{eqnarray*}
\calt \calt^{-1} & = & \calt^{-1} \calt = 1,\\
q^{2m} \calt^m c & = & c \calt^m,\\
\calt^m e_n & = & q^{-2(n+1)m} e_n \calt^m,\\
q^{2m} e_m c & = & c e_m,\\
q^{m-n} e_m e_n - q^{n-m} e_n e_m & = & [m-n] e_{m+n} + \delta_{m+n, 0}
\frac{[m-1][m][m+1]}{[2][3]\tangle{m}} c,
\end{eqnarray*}

\noindent where $\displaystyle [m] := \frac{q^m - q^{-m}}{q - q^{-1}}$
and $\tangle{m} := q^m + q^{-m}$ for all $m \in \Z$.

One can now compute the group of weights, as well as $\sg(h)$ for any
monomial word $h$ in the above alphabet. The following is proved using
the defining relations above.

\begin{prop}
Setup as above.
\begin{enumerate}
\item The group of weights is $\Gamma_{\mathcal{U}_q} = (R^\times,
\cdot)$ (so every finite subgroup $\Pi$ is cyclic). A weight $r \in
R^\times$ kills $c$ and all $e_n$, and sends $\calt$ to $r$.

\item $\sg(h) = 0$ if the monomial word $h$ contains $c$ or any $e_n$.
Moreover, $\sg(\calt^m) = 0$ unless $|\Pi|$ divides $m$, in which case
$\sg(\calt^m) = |\Pi|$.
\end{enumerate}
\end{prop}

\comment{
\begin{proof}\hfill
\begin{enumerate}
\item Given any weight $\mu$, evaluating on the second relation suggests
that $\mu(c) = 0$. Similarly, from the third relation, we get that
$\mu(e_n) = 0$ if $n \neq -1$. Now choose $m = 1, n = -2$, and evaluate
$\mu$ on the last relation. Then we get that $0 - 0 = [3] \mu(e_{-1})$,
so each $\mu \in \Gamma$ kills $c, e_n$ for all $n$. Finally, the first
relation says that $\mu(\calt) \in R^\times$.

Thus we now define $\mu_r \in \Gamma$, by sending $\calt$ to $r \in
R^\times$ and $c, e_n$ to zero. This is easily seen to be compatible with
all the defining relations, so the first part is proved.

\item The first statement is now clear. Next, if $\Pi$ is not the trivial
group, then $\Pi$ is nontrivial and cyclic, and $\calt^m \in \ggh$ if and
only if any generator $\la$ of $\Pi$ satisfies: $\la(\calt)^m = 1$. In
other words, the order of $\la$ (which equals $|\Pi|$) divides $m$. Now
invoke Proposition \ref{P5} to finish the proof.
\end{enumerate}
\end{proof}
}
\end{example}

\begin{example}[Quantum linear groups]
For the definitions, we refer to \cite[\S 2]{HaKr}. The {\it quantum
general} (respectively {\it special}) {\it linear group} $GL_q(n) =
R_q[GL_n]$ (respectively $SL_q(n) = R_q[SL_n]$) is the localization of
the algebra $\B$ (defined presently) at the central {\it quantum
determinant}
\[ \det{}_q := \sum_{\pi \in S_n} (-q)^{l(\pi)} \prod_{i=1}^n u_{i,
\pi(i)} \]

\noindent (respectively the quotient of $\B$ by the relation $\det_q =
1$). Here, the algebra $\B = R_q[\mathfrak{gl}(n)]$ is generated by $\{
u_{ij} : 1 \leq i,j \leq n \}$, with relations
\begin{eqnarray*}
&& u_{ik} u_{il} = q u_{il} u_{ik}, \qquad u_{ik} u_{jk} = q u_{jk}
u_{ik},\\
&& u_{il} u_{jk} = u_{jk} u_{il},\qquad u_{ik} u_{jl} - u_{jl} u_{ik} =
(q-q^{-1}) u_{il} u_{jk},
\end{eqnarray*}

\noindent where $q \in R^\times$, and $i < j,k < l$.

As above, it is possible to compute the group of weights $\Gamma$ for
both families of algebras, in the ``nontrivial" case $q \neq \pm 1$. In
either case, note that $\det_q \neq 0$. Given any permutation $\pi \in
S_n$, suppose there exist $i < j$ such that $\pi(i) > \pi(j)$. Then
\[ (q - q^{-1})\mu(u_{i, \pi(i)} u_{j, \pi(j)}) = \mu([u_{i, \pi(j)},
u_{j, \pi(i)}]) = 0, \]

\noindent whence $\mu(u_{i, \pi(i)} u_{j, \pi(j)}) = 0$ for all $\mu \in
\Gamma$.

The only permutations for which this does not happen is $\{ \pi \in S_n :
i<j \To \pi(i) < \pi(j) \} = \{ \id \}$. Hence $\mu(\det_q) = \prod_i
\mu(u_{ii}) \neq 0$, whence no $u_{ii}$ is killed by any $\mu$. But now
for $i<l$,
\[ \mu(u_{ii} u_{il}) = q \mu(u_{il} u_{ii}), \qquad \mu(u_{ii} u_{li}) =
q \mu(u_{li} u_{ii}), \]

\noindent whence $\mu(u_{ij}) = 0$ for all $i \neq j, \mu \in \Gamma$. In
particular, since $\dd(u_{ij}) = \sum_{k=1}^n u_{ik} \otimes u_{kj}$ for
(all $i,j$ and) both $GL_q(n)$ and $SL_q(n)$, hence $\Gamma_{GL} \cong
(R^\times)^n$, and $\Gamma_{SL} \cong (R^\times)^{n-1}$ (both under
coordinate-wise multiplication), since $\mu(u_{nn}) = \prod_{i=1}^{n-1}
\mu(u_{ii})^{-1}$ in the latter case.

Finally, computing $\sg$ now reduces to the above results and the first
example (of the free group $G = \mfgn{0}$). This is because any $h \in
GL_q(n)$ or $SL_q(n)$ can be reduced to a sum of monomial words, and such
a word is not killed by any $\mu$ if and only if there is no contribution
from any $u_{il}, i \neq l$.
\end{example}

\begin{example}[Hopf regular triangular algebras]
These were defined in \cite{Kh2} (in the special case $\Gamma = 1$).

\begin{defn}
An associative $k$-algebra $A$ (over a ground field $k$) is a {\em Hopf
RTA} (or {\em HRTA}), if:
\begin{enumerate}
\item The multiplication map : $B_- \otimes_k H \otimes_k B_+ \to A$ is
an isomorphism, for some (fixed) associative unital $k$-subalgebras $H,
B_\pm$ of $A$, and $H$ is, in addition, a commutative Hopf algebra.

\item The set $G := \Hom_{k-alg}(H,k)$ contains a free abelian group with
finite basis $\dd$, so that $B_\pm = \bigoplus_{\la \in \pm
\mathbb{Z}_{\geqslant 0} \dd} (B_\pm)_\la$. Each summand here is a
finite-dimensional weight space for the (usual) adjoint action of $H$,
and $(B_\pm)_0 = k$.

\item There exists an anti-involution $i$ of $A$, so that $i|_H =
\id|_H$.
\end{enumerate}
\end{defn}

\noindent This is a large family of algebras, that are widely studied in
representation theory. Examples (when $\ch(k) = 0$) are $U(\mfg)$ for
$\mfg$ a semisimple, symmetrizable Kac-Moody, (centerless Virasoro), or
(centerless) extended Heisenberg Lie algebra. Other examples include
quantum groups $U_q(\mfg)$ or (quantized) infinitesimal Hecke algebras
over $\mathfrak{sl}_2$.\medskip

In all these examples, the computations reduce to $H$:

\begin{lemma}
$\Gamma_A \subset \Gamma_H$.
\end{lemma}

The same happens if one works with {\em skew group rings} over $A$, for
example, wreath products $S_n \wr A := A^{\otimes n} \rtimes S_n$; in
this case $H$ is replaced in the above lemma, by the Hopf subalgebra
$H^{\otimes n} \rtimes S_n$ of $S_n \wr A$.

\begin{proof}
Use Proposition \ref{Pq}, with $V = N_- + N_+$, where $N_\pm$ are the
augmentation ideals in $B_\pm$. Hence $\Gamma_A \subset \Gamma_H$ (note
that each weight must also kill $[A,A] \cap H$).
\end{proof}
\end{example}

\begin{example}[Finite-dimensional pointed Hopf algebras]\label{Eas}
Assume that:
\begin{enumerate}
\item $R = \overline{k(R)}$ is an algebraically closed field of
characteristic zero.

\item $H$ is a finite-dimensional pointed Hopf algebra over $R$.

\item $G(H)$ is a (finite) abelian group of order coprime to 210.
\end{enumerate}

\noindent Then by the Classification Theorem 0.1 of \cite{AS}, $H$ is
generated by $G(H)$ and some skew-primitive (defined in the next section)
generators $\{ x_i \}$, that satisfy $g x_i g^{-1} = \chi_i(g) x_i$ for
all $i$ and all $g \in G(H)$. Since $\chi_i \neq \vi\ \forall i$ by
\cite[Equation (0.1)]{AS}, hence Proposition \ref{Pq} again reduces the
computations here, to the grouplike case: $\Gamma_H \subset
\Gamma_{G(H)}$.
\end{example}

\section{Skew-primitive elements}\label{S3}

As the previous example \ref{Eas} suggests, large families of Hopf
algebras that are the subject of much study in the literature, are
generated by grouplike and ``skew-primitive" elements. Hence in the rest
of this paper, we address the question of computing $\sg$ at (monomial
words in) such elements.

\begin{defn}
An element $h \in H$ is {\it skew-primitive} if $\dd(h) = g \otimes h + h
\otimes g'$ for grouplike $g, g' \in G(H)$. Denote the set of such
elements by $H_{g,g'}$. Then $g-g' \in H_{g,g'} \cap H_{g',g}$, and
$H_{1,1} = H_{prim}$ (recall Definition \ref{D1}).

An element $h$ is {\it pseudo-primitive} (respectively {\it almost
primitive}) with respect to $\Pi$ if, moreover, $g^{-1}g' \in \ggh$
(respectively $g,g' \in \ggh$). In future, we may not specify the finite
subgroup $\Pi$ of $\Gamma_H$, because it is part of the given data.
\end{defn}\smallskip

Thus,
$\{ \mbox{skew-primitive} \} \supset \{ \mbox{pseudo-primitive} \}
\supset \{ \mbox{almost primitive} \} \supset \{ \mbox{primitive} \}$.

\begin{remark}
Note that grouplike and skew-primitive elements in $H$ correspond to
$1$-dimensional and non-semisimple $2$-dimensional $H$-comodules,
respectively. To see this, suppose $V_0 = R v_0$ and $g \in H$. Define
$\rho(v_0) := v_0 \otimes g$. Then $\rho$ induces an $H$-comodule
structure on $V_0$ if and only if $g$ is grouplike.
Similarly, given $V_1 = R v_0 \oplus R v_1$ and $g,g',h \in H$, define:
\[ \rho(v_0) := v_0 \otimes g, \qquad \rho(v_1) = v_0 \otimes h
+ v_1 \otimes g'. \]

\noindent Then $\rho$ induces an $H$-comodule structure on $V_0$ if and
only if $g,g'$ are grouplike and $h \in H_{g,g'}$ is
$(g,g')$-skew-primitive.
\end{remark}

\begin{lemma}\label{L5}
Suppose $\Pi \subset \Gamma_H$ is a finite subgroup of weights, and $h
\in H_{g,g'}$ is as above.
\begin{enumerate}
\item The set $\{ \g \in \Gamma_H : \g(h) = 0 \}$ is a subgroup of
$\Gamma_H$.

\item $\vi(h) = 0$ and $S(h) = -g^{-1} h (g')^{-1} \in H_{(g')^{-1},
g^{-1}}$.

\item If $g_0$ is any grouplike element, then $g_0 h$ and $h g_0$ are
also skew-primitive. Moreover, $g_0 - g_0^{-1} \in H_{g_0, g_0^{-1}} \cap
H_{g_0^{-1}, g_0}$.

\item For all $n \geq 0$, one also has:
\begin{equation}\label{E1}
\dd^{(n)}(h) = \sum_{i=0}^n g^{\otimes i} \otimes h \otimes (g')^{\otimes
(n-i)}.
\end{equation}

\item 
For any $\g \in \Pi$, either $\g(g) \neq \g(g')$, or $\g(h) = 0$, or
$\ch(R) | n_\g$.

\item $h$ is pseudo-primitive if and only if $g^{-1}h, h g^{-1},
(g')^{-1}h, h (g')^{-1}$ are almost primitive. If $g_0 \in G(H)$ and $h$
is pseudo-primitive, then so are $g_0 h$ and $h g_0$.
\end{enumerate}
\end{lemma}

\begin{proof}
The first part is an easy verification. The second part follows from the
statements
\[ (\id * \vi)(h) := \mu(\id \otimes \vi)\dd(h) = h, \qquad (\id * S)(h)
:= \mu(\id \otimes S)\dd(h) = \vi(h). \]

\noindent The third, fourth and last parts are now easy to verify. For
the fifth part, suppose $\g(g) = \g(g')$ and $\g$ has order $n_\g \geq 1$
in $\Pi$. Then compute using equation \eqref{E1} above:
\[ 0 = \vi(h) = \g^{*n_\g}(h) = n_\g \g(h) \g(g)^{n_\g - 1}. \]

\noindent The result follows, since $\g(g)$ is a unit and $R$ is an
integral domain.
\end{proof}

Our main result here is to compute $\sg(g_1 h g_2)$ for any $g_1, g_2 \in
G(H)$, or equivalently by the lemma above, $\sg(h)$ for all
skew-primitive $h$.
Equation \eqref{E1} implies that if $\g \in \Gamma_H$, and $\g(g) \neq
\g(g')$, then
\begin{equation}\label{E2}
\g^n(h) = \g^{\otimes n}(\dd^{(n-1)}(h)) = \g(h) \cdot \frac{\g(g)^n -
\g(g')^n}{\g(g) - \g(g')}.
\end{equation}

\subsection{The main result}

In all that follows below, assume that $h$ is skew-primitive, with
$\dd(h) = g \otimes h + h \otimes g'$.

\begin{theorem}\label{T3}
A skew primitive $h \in H_{g,g'}$ satisfies at least one of the following
three conditions:
\begin{enumerate}
\item If there is $\la \in \Pi$ so that $\la(g), \la(g') \neq 1$, then
$\sg(h) = 0$. If no such $\la$ exists, then one of $g, g'$ is in $\ggh$.

\item Suppose only one of $g, g'$ is in $\ggh$, so that there exists $\la
\in \Pi$ with exactly one of $\la(g), \la(g')$ equal to 1. Then
$\displaystyle \sg(h) = \frac{|\Pi|\la(h)}{1 - \la(gg')}$.

\item If $\la(gg') = 1$ for all $\la \in \Pi$, then $\sg(h) = \sum_{\g
\neq \vi = \g^2} \g(h)$, and $2 \sg(h) = 0$.
\end{enumerate}
\end{theorem}

\begin{remark}\hfill
\begin{enumerate}
\item Thus, the expression $\la(h) / (1 - \la(gg')) = \sg(h) \in k(R)$ is
independent of $\la$ (as long as $\la(gg') \neq 1$), for such $h$. As the
proof indicates, $1 - \la(gg')$ should really be thought of as $1 -
\la(g)$ or $1 - \la(g')$ (depending on which of $g'$ and $g$ is in
$\ggh$).

Moreover, equation \eqref{E2} implies, whenever $\g^n(g) \neq \g^n(g')$,
that
\[ \frac{\g^n(h)}{\g^n(g) - \g^n(g')} = \frac{\g(h)}{\g(g) - \g(g')}. \]

\noindent It is not hard to show that both of these are manifestations of
the following easy fact:

\begin{lemma}\label{Llevel}
Given $h \in H_{g,g'}$, define $N_h := \{ \g \in H^* : \g(g) \neq \g(g')
\}$. Suppose $\mu, \la \in N_h$. Then $(\mu*\la)(h) = (\la*\mu)(h)$ if
and only if $f_h(\mu) = f_h(\la)$, where $f_h : N_h \to k(R)$ is given by
\[ f_h(\g) := \frac{\g(h)}{\g(g) - \g(g')}. \]
\end{lemma}

\noindent (In other words, weights commute at $h$ precisely when they lie
on the same ``level surface" for $f_h$.)

\item Also note that if the first two parts fail to hold, then both $g,
g' \in \ggh$, and the final part holds. Thus, the above theorem computes
$\sg(h)$ for all skew-primitive $h$, if $\ch(R) \neq 2$ or not both of
$g,g'$ are in $\ggh$. We address the case when $\ch(R) = 2$ and $g,g' \in
\ggh$, in the next subsection.
\end{enumerate}
\end{remark}

\begin{proof}[Proof of Theorem \ref{T3}]\hfill
\begin{enumerate}
\item Apply equation \eqref{E0} to $h$, to get
$\sg(h)
= \la(g) \sg(h) + \la(h) \sg(g')$.
Since $g' \notin \ggh$, hence the second term vanishes, and we are left
with $(1 - \la(g)) \sg(h) = 0$. But $\la(g) \neq 1$.

\comment{
By Lemma \ref{L5} above, there are three cases:
\begin{enumerate}
\item If $\la(h) = 0$, then we are done by Proposition \ref{P6} above.

\item If $\la(g) \neq \la(g')$, then we use equation \eqref{E4}, and
since the last term on the right side vanishes, we are again done.

\item Finally, if $\la(g) = \la(g')$ and $\la(h) \neq 0$, then $n = n_\la
= |\tangle{\la}|$ is divisible by $\ch(R)$. Now use equation \eqref{E1}
to compute:
\[ \sla(h) = \sum_{i=1}^n i \la(h) \la(g)^{i-1} \]

\noindent Since $(a-1)^2 \sum_{i=1}^m i a^{i-1} = m a^{m+1} - (m+1) a^m +
1$ for any $m \in \N, a \in R$, and using that $n=0$ in $R$, we get
\[ \sla(h) = \frac{1 - \la(g)^n}{(1 - \la(g))^2} \]

\noindent and since $\la^n = \vi$, we get that $\sla(h) = 0$. Now choose
a set $\B$ of coset representatives as above, and compute:
\[ \sg(h) = \sum_{\beta \in \B} \left( \beta(g) \cdot \sla(h) + \beta(h)
\cdot \sla(g') \right) \]

\noindent and the second factor in either summand vanishes, using
Proposition \ref{P5} above.\\
\end{enumerate}
}

Next, if no such $\la \in \Pi$ exists, then $\Pi = \Pi_g \cup \Pi_{g'}$,
where $\Pi_g, \Pi_{g'}$ were defined before Proposition \ref{P5}. We
claim that one of the two sets is contained in the other, whence one of
$g, g'$ is in $\ggh$. For if not, then choose $\g, \g' \in \Pi$ so that
neither $\g(g)$ nor $\g'(g')$ equals 1 (whence $\g(g') = 1 = \g'(g)$).
Then one verifies that $\g \g' \notin \Pi_g \cup \Pi_{g'}$.
and this is a contradiction. Thus one of $\Pi_g \setminus
\Pi_{g'},\ \Pi_{g'} \setminus \Pi_g$ is empty.\medskip

\item Suppose $g' \in \ggh, \la(g) \neq 1$ for some $\la$ (the other case
is similar). Now apply equation \eqref{E0} and Proposition \ref{P5}, and
compute:
\[ \sg(h) = \la(h) \frac{\sg(g) - \sg(g')}{\la(g) - \la(g')} = \la(h)
\frac{0 - |\Pi|}{\la(g) \cdot 1 - 1} = \frac{\la(h) |\Pi|}{1 - \la(g)
\la(g')}. \]

\item If $\la(gg') = 1$ for any $\la \in \Pi$, then: $\la^{-1}(h) =
-\la(h)$.
%
Thus $\la + \la^{-1}$ kills $h$, and the first equation now follows
because $\vi(h) = 0$. The second is also easy:
$2 \sg(h) = \sum_{\g \in \Pi}(\g(h) + \g^{-1}(h)) = 0$.
\end{enumerate}
\end{proof}

\comment{
For the second part, we have a special $\la$, and partitioning
$\tangle{\la}$ into subsets as above, $\sla(h) = 0$. As above, we now
choose a set $\B$ of coset representatives for $\tangle{\la}$ in
$\Gamma_H$, and compute:
\[ \sg(h) = \sum_{\beta \in \B,\ \g \in \tangle{\la}} (\beta \otimes
\g)(\dd(h)) = \sum_{\beta \in \B} \left( \beta(g) \cdot \sla(h) +
\beta(h) \cdot \sla(g') \right) \]

\noindent and this vanishes because the second factor in each term of the
summand is zero (by above and by the previous subsection).\\
}

\subsection{The characteristic 2 case}

The only case that Theorem \ref{T3} does not address, is when $\ch(R) =
2$ and $g,g' \in \ggh$. We now address this case.

\begin{prop}\label{Punproved}
Suppose $\Pi$ is as above, $\ch(R) = 2$, and $h \in H_{g,g'}$ is almost
primitive with respect to $\Pi$.
\begin{enumerate}
\item If $\Pi$ has odd order, then $\sg(h) = 0$.

\item If 4 divides $|\Pi|$, then $\sg(h) = 0$.

\item If $\Pi$ has even order but $4 \nmid |\Pi|$, then $\sg(h) = \g(h)$
for any $\g \in \Pi$ of order exactly 2. This may assume any nonzero
value in $R$.
\end{enumerate}
\end{prop}

\noindent We omit the proof, since this result is a special case of more
general results in general (positive) characteristic, which we state and
prove later. See Theorems \ref{T6} and \ref{T8}, as well as Remark
\ref{R711}.

\subsection{A degenerate example: quiver (co)algebras}

We conclude this section with an example that is not a Hopf algebra, but
an algebra with coproduct, and is generated by grouplike and
skew-primitive elements. (Thus, $\Pi$ is no longer necessarily a group of
weights, but a semigroup.)

Consider a quiver $Q = (Q_0, Q_1)$, which is a directed graph with vertex
set $Q_0$ and edges $Q_1$. Thus, there exist source and target maps $s,t
: Q_1 \to Q_0$ such that every edge $e \in Q_1$ starts at $s(e)$ and ends
at $t(e)$. A path in $Q$ is a finite sequence of edges $a_1 \cdots a_n$
such that $t(a_i) = s(a_{i+1})$. We also write $s(p) := s(a_1)$ and $t(p)
:= t(a_n)$, and the length of the path is said to be $n$. Vertices $v \in
Q_0$ are paths of length zero, and one writes $s(v) = t(v) := v$.

There are two structures on the free $R$-module $R Q$ with basis
consisting of all paths in $Q$. The {\it path algebra} is defined by
setting the product of two paths $a_1 \cdots a_n$ and $b_1 \cdots b_m$ to
be their concatenation $a_1 \cdots a_n b_1 \cdots b_m$ if $t(a_n) =
s(b_m)$; otherwise the product is zero. Then $R Q$ is an associative
$R$-algebra which contains enough idempotents $\{ v : v \in Q_0 \}$.
However, $R Q$ contains a unit if and only if $Q_0$ is finite, in which
case the unit is $\sum_{v \in Q_0} v$.

Another structure on $R Q$ is given by defining $\Delta : R Q \to R Q
\otimes R Q$ via: $\Delta(p) := \sum_{(q,r) : qr = p} q \otimes r$. Also
define $\varepsilon : R Q \to R$ via: $\varepsilon(v) := 1$ for all $v
\in Q_0$ and $\varepsilon(p) := 0$ for all paths $p$ of positive length.
This structure makes $R Q$ into the {\it path coalgebra}.

The path (co)algebra has been the subject of much study in the
literature, and it is natural to ask for which quivers $Q$ are these two
structures on $R Q$ compatible. It is not hard to show that the answer
is: very few.

\begin{lemma}
The coproduct $\Delta$ is multiplicative if and only if there are no
paths of length $ \geq 2$.
\end{lemma}

\noindent Thus, the quiver bialgebra is a ``degenerate" example. Now
suppose $Q$ has no paths of length $\geq 2$ (and hence, no self-loops).
Then $\Delta$ is indeed multiplicative. Nevertheless, if $1 < |Q_0| <
\infty$, then the unit in $R Q$ is not grouplike. It is now possible to
show the following.

\begin{prop}
There exists a bijection from $Q_0 \coprod \{ 0 \}$ to the semigroup
$\Gamma_{R Q}$, sending $v$ to $\lambda_v$ that sends $v$ to $1$ and all
other paths to $0$.
\end{prop}

\begin{proof}
Given $v \neq v' \in Q_0$ and $\lambda \in \Gamma_{R Q}$,
\[ \lambda(v) = \lambda(v^2) = \lambda(v)^2, \qquad \lambda(v)
\lambda(v') = \lambda(vv') = 0. \]

\noindent Since $R$ is a unital integral domain, this implies that at
most one $\lambda(v)$ is nonzero - and then it equals $1$. Moreover, if
$p$ is any path of length $1$ in $Q_1$, then at least one of its vertices
is killed by $\lambda$, whence so is $p$. This supplies the desired
bijection. In particular, $\Gamma_{R Q} \cong Q_0 \coprod \{ 0 \}$ as
semigroups, with composition given by:
\[ \lambda_v * \lambda_{v'} = \delta_{v,v'} \lambda_v, \qquad \lambda_v *
0 = 0 * \lambda_{v'} = 0, \qquad \forall v,v' \in Q_0. \]
\end{proof}

In fact, the path algebra and path coalgebra are dual to one another,
morally speaking. More precisely, each of them can be ``recovered" inside
the dual space of the other; see \cite{DIN} for more details.

\section{Subgroups and subquotients of groups of weights}\label{Ssubquot}

\subsection{Subgroups associated to arbitrary elements}

Recall that the goal of this article is to compute $\sg$ at any element
$h$ in a Hopf algebra $H$. We start this section with the following
constructions.

\begin{defn}
Suppose $R$ is a commutative unital integral domain.
\begin{enumerate}
\item Suppose $H$ is an $R$-algebra, such that $\Gamma = \Gamma_H$ has a
group structure $*$ on it. Define $\Gamma_h \subset \Gamma$ to be the
subgroup of $\Gamma$ that ``stabilizes" $h$. In other words,
\[ \Gamma_h := \{ \g \in \Gamma : (\beta * \g * \delta)(h) = (\beta *
\delta)(h)\ \forall \beta, \delta \in \Gamma \}. \]

\item Given a coalgebra $H$, and $h \in H$, define $C_h$ to be the
$R$-subcoalgebra generated by $h$ in $H$.

\item Given a Hopf algebra $H$, define $\Gamma'_h$ to be the {\it fixed
weight monoid} of $h$, given by
$\Gamma'_h := \{ \g \in \Gamma : \g|_{C_h} = \vi|_{C_h} \}$.
\end{enumerate}
\end{defn}

\noindent In particular, $\g(h) = \vi(h)$ if $\g \in \Gamma_h$.\medskip

A later subsection will discuss how this allows us to consider
subquotients of $\Gamma$; but first, here are some observations involving
these subgroups.

\begin{prop}\label{Pstab}
Suppose $H$ is a Hopf algebra, and $h \in H$.
\begin{enumerate}
\item For all $h,\ \Gamma_h$ is a normal subgroup of $\Gamma$, and
$\Gamma'_h \subset \Gamma_h$ is a monoid closed under
$\Gamma$-conjugation.

\item Given $\{ h_i : i \in I \} \subset H$, and $h \in \tangle{h_i}$
(i.e., in the subalgebra generated by the $h_i$'s), $\Gamma_h \supset
\bigcap_{i \in I} \Gamma_{h_i}$, and similarly for the $\Gamma'$s.

\item Given any $h_i \in H$ (finitely many), suppose $\Pi = \times_i
\Pi_i$, with $\Pi_i \subset \Gamma'_{h_j}$ whenever $i \neq j$. Then
$\displaystyle \sg(\bfh) = \prod_i \Sigma_{\Pi_i}(h_i)$.
\end{enumerate}
\end{prop}

\begin{proof}
The first and third parts are straightforward computations.
\comment{
We show that $\Gamma_h$ is closed under multiplication, taking inverses,
and $\Gamma$-conjugation respectively; for convenience, we omit the
$'*'$s and the $(h)$. Thus, as far as evaluating at $h$ goes, if $\g,\g'
\in \Gamma_h$, and $\beta, \delta, \mu \in \Gamma$, then
\begin{eqnarray*}
\beta \g \g' \delta & = & \beta \g (\g' \delta) = \beta \g' \delta =
\beta \delta\\
\beta \cdot \g^{-1} \delta & = & \beta \g^2 (\g^{-1} \delta) = \beta \g
\delta = \beta \delta\\
\beta (\mu \g \mu^{-1}) \delta & = & (\beta \mu) \g (\mu^{-1} \delta) =
\beta \mu \cdot \mu^{-1} \delta = \beta \delta
\end{eqnarray*}

Next, we of course have $\vi \in \Gamma'_H$, and $\Gamma'_h \subset
\Gamma_h$ since given $\g' \in \Gamma'_h$ and $\beta, \delta \in \Gamma$,
we have
\[ (\beta * \g' * \delta)(h) = \sum \beta(\one{h}) \g'(\two{h})
\delta(\three{h}) = \sum \beta(\one{h}) \vi(\two{h}) \delta(\three{h}) \]

\noindent and this equals $(\beta * \delta)(h)$. Finally, $\Gamma'_h$ is
closed under multiplication and $\Gamma$-conjugation, because given
$\beta \in \Gamma,\ \g',\g'' \in \Gamma'_h$, and $h' \in C_h$, we get
that $(\beta * \g' * \beta^{-1})(h')$ equals
\[ \sum \beta(\one{h'}) \g'(\two{h'}) \beta^{-1}(\three{h'}) = \sum
\beta(\one{h'}) \vi(\two{h'}) \beta(S(\three{h'})) = \vi(h') \]

\noindent and
\[ (\g' * \g'')(h) = \sum \g'(\one{h}) \g''(\two{h}) = \sum \vi(\one{h})
\vi(\two{h}) = \vi(h) \]
}
For the second part, for all $\beta, \delta \in \Gamma,\ \g \in \bigcap_i
\Gamma_{h_i}$, and polynomials $p$ in the $h_i$'s,
\[ (\beta * \g * \delta)(p(h_i)) = p((\beta * \g * \delta)(h_i)) =
p((\beta * \delta)(h_i)) = (\beta * \delta)(p(h_i)). \]

\noindent The outer equalities hold because weights are algebra maps.

The proof for the $\Gamma'$s is as follows: if $h = p(h_i)$ as above,
then since $\dd$ is multiplicative, hence any $h' \in C_h$ is expressible
as a polynomial in elements $h'_j \in \cup_i C_{h_i}$ - say $h' =
q(h'_j)$. In particular, if $\g \in \Gamma'_{h_i}$ for all $i$, then
\[ \g(q(h'_j)) = q(\g(h'_j)) = q(\vi(h'_j)) = \vi(q(h'_j)) \]

\noindent where once again, the outer equalities hold because weights are
algebra maps. In other words, $\g(h') = \vi(h')$.
\comment{
\item Given $\g_i \in \Pi_i\ \forall i$, we compute (for each $j$):
\begin{eqnarray*}
(*_i \g_i)(h_j) & = & \sum \g_1(\one{(h_j)}) \g_2(\two{(h_j)}) \dots
\g_n((h_j)_{(n)})\\
& = & \sum \g_j((h_j)_{(j)}) \prod_{i \neq j} \vi((h_j)_{(i)}) =
\g_j(h_j)
\end{eqnarray*}

\noindent so that we have
\[ \sg(h) = \sum_{(\g_i)_i \in \Pi} (*_i \g_i)(\bfh) = \sum_{(\g_i)_i \in
\Pi} \prod_i \g_i(h_i) = \prod_{i=1}^n \sum_{\g_i \in \Pi_i} \g_i(h_i) =
\prod_{i=1}^n \Sigma_{\Pi_i}(h_i) \]
}
\end{proof}

We also mention two examples; the proofs are straightforward.

\begin{lemma}\label{Lstab}
If $g \in G(H)$, then this definition of $\Gamma_g$ coincides with the
previous one: $\Gamma_g = \Gamma'_g = \{ \g \in \Gamma : \g(g) = 1 \}$.
If $h \in H_{g,g'}$, then $\Gamma_h \subset \Gamma_g \cap \Gamma_{g'}$,
or $\Gamma_h = \Gamma$. In both cases, $\Gamma'_h = \Gamma_h \cap
\Gamma_g \cap \Gamma_{g'}$.
\end{lemma}

\comment{
\begin{proof}
The first statement is easy to show. Now if $h \in H_{g,g'}$ and $\g \in
\Gamma_h$, then $\g(h) = \vi(h) = 0$, whence the definition of $\Gamma_h$
also gives us that
\[ \beta(h) \g(g) = \beta(h) = \beta(h) \g(g')\ \forall \beta \in
\Gamma,\ \g \in \Gamma_h. \]

\noindent Therefore either $\beta(h) = 0$ for all $\beta$ - whence we get
that $\Gamma_h = \Gamma$ - or some $\beta(h) \neq 0$, whence $\g(g) =
\g(g') = 1$ for all $\g \in \Gamma_h$, as claimed. It is now easy (in
both cases) to compute $\Gamma'_h$.
\end{proof}\medskip
}

\subsection{Subquotients}

We now compute $\sg(h)$ for more general $\Pi$. The following result is
used later.

\begin{lemma}\label{L1}
Fix $h \in H$, and choose any subgroup $\Gamma'$ of $\Gamma_h$ that is
normal in $\Gamma$. Also fix a finite subgroup $\Pi$ of $\Gamma
/\Gamma'$. Now fix any lift $\widetilde{\Pi}$ of $\Pi$ to $\Gamma$, and
define
\[ \sg(h) := \sum_{\g'' \in \widetilde{\Pi}} \g''(h), \qquad \Pi^\circ :=
\{ \g \in \Gamma : (\g + \Gamma') \in \Pi \subset \Gamma / \Gamma' \}. \]
\begin{enumerate}
\item $\sg(h)$ is well-defined.
\item If a subgroup $\Gamma'' \subset \Gamma_h$ is normal in $\Pi^\circ$
(e.g., $\Gamma'' = \Gamma',\ \Gamma_h \cap \Pi^\circ$), then
\[ \Sigma_{\Pi^\circ}(h) = |\Gamma''| \Sigma_{\Pi^\circ / \Gamma''}(h). \]
\end{enumerate}
\end{lemma}

\begin{proof}
For the first part, choose any other lift $\Pi'$ of $\Pi$. If $\g' \in
\Pi', \g'' \in \widetilde{\Pi}$ are lifts of $\g \in \Pi$, then
$(\g')^{-1} * \g'' \in \Gamma_h$, so
\[ \sum_{\g'' \in \widetilde{\Pi}} \g''(h) = \sum_{\g' \in \Pi'} (\g' *
((\g')^{-1} * \g''))(h) = \sum_{\g' \in \Pi'} \g'(h) \]

\noindent by definition of $\Gamma_h$. The other part is also easy, since
by the first part, $\Sigma_{\Pi^*}(h)$ is also well-defined, where the
finite group is $\Pi^* := \Pi^\circ / \Gamma''$. If $\Pi^{**}$ is any
lift to $\Pi^\circ$ of $\Pi^*$, then
\[ \Sigma_{\Pi^\circ}(h) = \sum_{\g' \in \Pi^{**},\ \beta \in \Gamma''}
(\g' * \beta)(h) = \sum_{\g' \in \Pi^{**},\ \beta \in \Gamma''} \g'(h) =
|\Gamma''| \sum_{\g' \in \Pi^{**}} \g'(h), \]

\noindent and this equals the desired amount.
\end{proof}

While a special case of the second part is that $\Sigma_{\Pi^\circ}(h) =
|\Gamma'| \sg(h)$, we really use the result when $\Gamma'$ is itself
finite, and we replace $\Pi^\circ$ by $\Pi$. The equation is then used to
compute $\sg(h)$.

\subsection{Pseudo-primitive elements}

For the rest of this paper, $H$ is an $R$-Hopf algebra, unless stated
otherwise. If $h$ is grouplike or (pseudo-)primitive, then it is easy to
see that $\sg(h) = |[\Pi, \Pi]| \Sigma_{\Pi_{ab}}(h)$, where $\Pi_{ab} :=
\Pi / [\Pi,\Pi]$ (we show the pseudo-primitive case presently). Thus, in
the grouplike case, the question of whether or not $h \in \ggh$ reduces
to evaluating (any lift of a set of) generators of the finite abelian
group $\Pi_{ab}$, at $h$.

\begin{prop}\label{P8}
Suppose $h \in H_{g,g'}$ is pseudo-primitive with respect to $\Pi$.
\begin{enumerate}
\item Then $(\g * \nu)(h) = (\nu * \g)(h) = \nu(g)\g(h) + \g(g)\nu(h)$
for all $\g, \nu \in \Pi$.


\item $\Gamma'_h \supset [\Pi,\Pi]$.

\item For any $m \geq 0$, $\g^{*m}(h) = m \g(g)^{m-1} \g(h)$ if $\g \in
\Pi$. In particular, if $\ch(R) = p$ is prime, then $\g^{*p}(h) = 0\
\forall \g \in \Pi$.
\end{enumerate}
\end{prop}

\begin{proof}
The first and last parts are by definition and induction respectively. As
for the second part, one shows the following computation for any
skew-primitive $h \in H_{g,g'}$, and $\beta, \beta' \in \Gamma_H$:
\begin{eqnarray*}
(\beta * \beta' * \beta^{-1} * (\beta')^{-1})(h) & = & \beta(h
(g')^{-1})(1 - \beta'(g (g')^{-1}))\\
&& + \beta'(h(g')^{-1})(\beta(g (g')^{-1}) - 1)
\end{eqnarray*}

\noindent using Lemma \ref{L5}. Since $h$ is pseudo-primitive with
respect to $\Pi$, this shows that every generator (and hence element)
$\la$ of $[\Pi,\Pi]$ satisfies: $\la(g) = \la(g') = 1$ and $\la(h) = 0$.
Since $C_h = Rh + Rg + Rg'$, these imply that $\la \in \Gamma'_h$.
\end{proof}

\comment{
\begin{enumerate}
\item This holds because $\nu(g) = \nu(g')$ for all $\nu \in \Pi$ (since
$h$ is pseudo-primitive).

\item This is because we note from Lemma \ref{L5} that
\[ \dd^{(3)}(h) = h \otimes (g')^{\otimes 3} + g \otimes h \otimes
(g')^{\otimes 2} + g^{\otimes 2} \otimes h \otimes g' + g^{\otimes 3}
\otimes h \]

\noindent We evaluate the expression $\beta * \beta' * \beta^{-1} *
(\beta')^{-1}$ at this, using Lemma \ref{L5} again. At times, we write
$\g(g) = \g(g')$ for $\g = \beta^{\pm 1}, (\beta')^{\pm 1}$, since $h$ is
pseudo-primitive. We thus get
\begin{eqnarray*}
& = & \beta(h) \beta(g')^{-1} + \beta(g) \beta'(h) \beta(g')^{-1}
\beta'(g')^{-1}\\
& - & \beta(g) \beta'(g) \beta(g^{-1} h (g')^{-1}) \beta'(g')^{-1}
- \beta(g) \beta'(g) \beta(g')^{-1} \beta'(g^{-1} h (g')^{-1})\\
& = & \beta(h) \beta(g)^{-1} + \beta'(h) \beta'(g)^{-1} - \beta(h)
\beta(g)^{-1} - \beta'(h) \beta'(g)^{-1}
\end{eqnarray*}

\item This is easy to show by induction.
\end{enumerate}
}

An easy consequence of Propositions \ref{Pstab} and \ref{P8} is

\begin{cor}\label{C1}
If $h \in H$ is (in the subalgebra) generated by grouplike and
pseudo-primitive elements (with respect to $\Pi$), then $\Gamma'_h
\supset [\Pi,\Pi]$.
\end{cor}

Also note that given some $h \in H$, one can compute $\sg(h)$ for more
general $\Pi$, and hence the results in this paper can be generalized;
however, we stay with the original setup when $\Pi \subset \Gamma$ (i.e.,
$\Gamma' = \{ \vi \}$). This (general) case is noteworthy, however,
because it is used below.\medskip

We conclude by specifying more precisely, what is meant by $\sg(h_1 \dots
h_n)$ for ``pseudo-primitive" $h_i$'s, when $\Pi$ is a subquotient of
$\Gamma$ as above. In this case, start with some skew-primitive $h_i$'s,
then let $\Pi$ be a finite subgroup of $\Gamma / \Gamma'$, for some
subgroup $\Gamma' \subset \cap_i \Gamma_{h_i}$ that is normal in
$\Gamma$. Moreover, if $\dd(h_i) = g_i \otimes h_i + h_i \otimes g'_i$,
then assume further that $\g(g_i) = \g(g'_i)$, for all elements $\g$ of
the subgroup $\Pi^\circ$ (defined above).

This is what is meant in the case of general $\Pi$, when we say that
$g_i, g'_i \in \ggh$ - i.e., that the $h_i$'s are pseudo-primitive (with
respect to $\Pi$). Similarly, to say that the $h_i$'s are almost
primitive with respect to $\Pi$ means that $\g(g_i) = \g(g'_i) = 1\
\forall \g \in \Pi^\circ$.

\section{Products of skew-primitive elements}\label{S5}

We now mention some results on (finite) products of skew-primitive
elements and grouplike elements. From now on, $\Pi$ denotes a finite
subgroup of $\Gamma = \Gamma_H$ and not a general subquotient; however,
in the next section, we need to use a subquotient $\Phi$ of this $\Pi$.

Since the set of skew-primitive elements is closed under multiplication
by grouplike elements, any ``monomial" in them can be expressed in the
form $\bfh = \prod_i h_i$. The related ``grouplike" elements that would
figure in the computations, are $\bfg = \prod_i g_i$ and $\bfg' = \prod_i
g'_i$.

\begin{stand}\label{D51}
For this section and the next two, assume that $h_i \in H_{g_i, g'_i}$
for all (finitely many) $i$.
\end{stand}\medskip

First, here are some results that hold in general.

\begin{prop}
If $\la(h_i) = 0\ \forall i$ for some $\la \in \Pi$, and (at least) one
of $\la(\bfg), \la(\bfg')$ is not 1, then $\sg(\bfh) = 0$.
\end{prop}

\begin{proof}
If $\la(h_i) = 0$ then $\la^m(h_i) = 0$ for all $i,m$. Now choose a set
$\B$ of coset representatives for $\tangle{\la}$ in $\Pi$, and assume
that $\la(\bfg') \neq 1$ (the other case is similar). Then compute:
\[ \sg(\bfh) = \sum_{\beta \in \B,\ \g \in \tangle{\la}} \prod_i
\left( \beta(g_i) \g(h_i) + \beta(h_i) \g(g'_i) \right). \]

\noindent Since $\g(h_i) = 0$ for all $\g \in \tangle{\la}$ and all $i$,
hence the entire product in the summand collapses, to give
$\sum_{\beta \in \B} \beta(\bfh) \cdot \sla(\bfg')$. But now
the second factor vanishes by our assumption (and Proposition \ref{P5}).
\end{proof}\medskip

Next, if $\sg(\prod_{i=1}^n h_i)$ is known for all skew-primitive
$h_i$'s, then one can evaluate the product of $(n+1)$ such $h$'s in some
cases. The following result relates $\sg$-values of strings to the
$\sg$-values of proper substrings (with skew-primitive ``letters"), that
are ``corrected" by grouplike elements. The proof is that both equations
below follow by evaluating equation \eqref{E0} at $\bfh$.

\begin{prop}
Suppose one of $\bfg, \bfg'$ is not in $\ggh$. Thus, if $\la(\bfg) \neq
1$ (or respectively $\la(\bfg') \neq 1$) for some $\la \in \Pi$, then
\[ \sg(\bfh) = \frac{ \sum_{\nu \in \Pi} \prod_i (\la(g_i) \nu(h_i) +
\la(h_i) \nu(g'_i)) - \la(\bfg) \sg(\bfh)}{1 - \la(\bfg)}, \]

\noindent or respectively,
\[ \sg(\bfh) = \frac{ \sum_{\nu \in \Pi} \prod_i (\nu(g_i) \la(h_i) +
\nu(h_i) \la(g'_i)) - \sg(\bfh) \la(\bfg')}{1 - \la(\bfg')}. \]
\end{prop}

Note here that both numerators on the right side have an $\sg(\bfh)$ in
them, which cancels the only such term present in the summations. Thus,
what one is left with in either case, are linear combinations of
$\sg$-values of ``corrected" proper substrings, with coefficients of the
form $\la(\prod g_j h_k)$.

Also note that if the $\sg$-values of all ``corrected" proper substrings
are known, and $\ch(R) \neq 2$, then the two propositions, one above and
one below, can be used, for instance,
to compute $\sg$ at all monomials of odd length (in skew-primitive
elements).\medskip

The statement and proof of the following result are essentially the same
as those of the last part of Theorem \ref{T3} above.

\begin{prop}
Suppose $\bfg \bfg' \in \ggh$, and $\ch(R) = 2$ if the number of $h_i$'s
is even. Then $\sg(\bfh) = \sum_{\g \neq \vi = \g^2} \g(\bfh)$, and $2
\sg(\bfh) = 0$.
\end{prop}

\noindent This is because once again, one shows that $\la^{-1}(\bfh) =
-\la(\bfh)\ \forall \la$.\medskip

The next result in this subsection is true for almost all values of $\ch
R$. The proof is immediate from the penultimate part of Lemma \ref{L5}
above.

\begin{prop}\label{P7}
If $\ch(R) \nmid |\Pi|$ and $g_i^{-1} g'_i \in \ggh$ (i.e., $h_i$ is
pseudo-primitive) for some $i$, then $\g(h_i) = 0\ \forall \g \in \Pi$.
In particular, $\sg(\bfh) = 0$.
\end{prop}

We conclude this section with one last result - in characteristic $p$.

\begin{theorem}\label{Tskew}
Suppose $\ch(R) = p > 0$, and $h_i \in H_{g_i, g'_i}$ for all $i$. Choose
and fix a $p$-Sylow subgroup $\Pi_p$ of $\Pi$.
\begin{enumerate}
\item Then each $h_i$ is almost primitive with respect to $\Pi_p$.

\item If $\Pi_p$ contains an element of order $p^2$, then $\sg(\bfh) =
\Sigma_{\Pi_p}(\bfh) = 0$.
\end{enumerate}
\end{theorem}

\noindent It is also shown later that $\sgp(\bfh) = 0$ whenever $\Pi_p
\ncong (\Z / p\Z)^k$ for any $k > 0$.

\begin{proof}\hfill
\begin{enumerate}
\item If $z^p = 1$ in $R$, then $(1-z)^p =
1 - z^p \mod p = 0$ in $R$. Since $R$ is an integral domain,
$z=1$.
Now assume that $|\Pi_p| = p^f$. Then $\g^{* p^f} = \vi$ for each $\g \in
\Pi_p$, whence $\g(g_i)^{p^f} = \g(g'_i)^{p^f} = 1$ for all $i$.
Successively set $z = \g(g_i)^{p^t}$, for $t = f-1, f-2, \dots, 1, 0$.
Hence $\g(g_i) = 1$ for all $i$; the other case is the same.\medskip

\item Next, if $\la$ has order $p^2$, then choose a set $\B$ of coset
representatives for $\tangle{\la}$ in $\Pi$, and compute using
Proposition \ref{P8} above (since all $h_i$'s are pseudo/almost primitive
with respect to $\Pi_p$, hence for $\tangle{\la}$):
\[ \sg(\bfh) = \sum_{\beta \in \B} \sum_{j=0}^{p^2-1} \prod_{i=1}^n
(\beta * \la^{*j})(h_i) = \sum_{\beta \in \B} \sum_{j=0}^{p^2-1}
\prod_{i=1}^n \left( \beta(h_i) \la^j(g_i) + j \beta(g_i) \la(g_i)^{j-1}
\la(h_i) \right). \]

\noindent Call the factor in the product $a_{i,j}$ (it really is
$a_{\beta, i, j}$). Now observe that $a_{i, p+j} = \la^p(g_i) a_{i,j}$
for all $i,j$, whence $a_{i, kp+j} = \la^{kp}(g_i) a_{i,j}$. Therefore
for any $\beta \in \B$, one can take $\sum_{k=0}^{p-1} \la(\bfg)^{kp}$
out of the summand. But $\la(\bfg)^p = 1$, so every $\beta$-summand
vanishes.
\comment{
in the above equation is
\[ \sum_{j=0}^{p^2 - 1} \prod_{i=1}^n a_{i,j} = \sum_{j=0}^{p - 1}
\sum_{k=0}^{p-1} \prod_{i=1}^n a_{i, kp+j} = \sum_{j=0}^{p-1}
\sum_{k=0}^{p-1} \la(\bfg)^{kp} \prod_{i=1}^n a_{ij}. \]

\noindent Finally, observe that the sum over $k$ yields $\sum_{k=0}^{p-1}
\ell^k$, where $\ell = \la(\bfg)^p = 1$ by the previous part. Thus the
sum is $p=0$, whence every $\beta$-summand above vanishes, and we are
done.
}
\end{enumerate}
\end{proof}

\section{Special case - abelian group of weights}\label{S6}

In this section, the focus is on evaluating $\sg(\bfh)$ in the special
case where the group $\Pi$ of weights is abelian.

\begin{defn}\hfill
\begin{enumerate}
\item For all $n \in \N$, define $[n] := \{ 1, 2, \dots, n \}$.

\item Given $I \subset [n]$, define $g_I := \prod_{i \in I} g_i$, and
similarly define $g'_I, h_I$.

\item Define $\Pi_p$ to be any fixed $p$-Sylow subgroup of $\Pi$ if
$\ch(R) = p>0$, and $\{ \vi \}$ otherwise. Also choose and fix a
``complementary" subgroup $\Pi'$ to $\Pi_p$ in $\Pi$ (if $\Pi$ is
abelian), i.e., $|\Pi_p| \cdot |\Pi'| = |\Pi|$. (And if $\ch(R) = 0$, set
$\Pi' := \Pi$.)
\end{enumerate}
\end{defn}

We now present two results. The first is (nontrivial only) when $\ch(R)$
divides the order of $\Pi$, and the second (which really is the main
result) is when it does not.

\begin{theorem}
Suppose $\Pi$ is abelian; let $\Pi_p, \Pi'$ be as above. Let $J \subset
[n]$ be the set of $i$'s such that $h_i$ is pseudo-primitive with respect
to $\Pi$. Then
\[ \sg(\bfh) = \Sigma_{\Pi'}(g_J h_{\notj}) \cdot \sgp(h_J). \]
\end{theorem}

\noindent In particular, if $\ch(R) = 0$, or $0< \ch(R) \nmid |\Pi|$,
then $\sgp(h_J) = \vi(h_J) = 0$, whence $\sg(\bfh) = 0$ too.

\begin{proof}
First note that $\la(g_i) = \la(g'_i) = 1$ for all $\la \in \Pi_p$ and
all $i$, by Theorem \ref{Tskew} above.
Since $\Pi$ is abelian, every $\g \in \Pi$ is uniquely expressible as $\g
= \beta * \la$ with $\beta \in \Pi',\ \la \in \Pi_p$. We now compute
$\g(h_i)$ in both cases: $i \in J$ and $i \notin J$.

First consider the case when $i \in J$. Then $\beta(h_i) = 0$ for all
$\beta \in \Pi'$, by Proposition \ref{P7}. Thus,
$\g(h_i) = \beta(g_i) \la(h_i) + \beta(h_i) \la(g_i) = \beta(g_i)
\la(h_i)$.

Now suppose $i \notin J$. Choose $\g \in \Pi$ so that $\g(g_i) \neq
\g(g'_i)$. Then
$(\g * \la)(h_i) = (\la * \g)(h_i)$,
which leads (upon simplifying) to
\[ \la(h_i) (\g(g_i) - \g(g'_i)) = \g(h_i)(\la(g_i) - \la(g'_i)). \]

\noindent But $\la(g_i) = \la(g'_i)$, and $\g(g_i) \neq \g(g'_i)$, so
$\la(h_i) = 0$ for all $\la \in \Pi_p$. Hence:
\[ \g(h_i) = (\beta * \la)(h_i) = \beta(g_i) \la(h_i) + \beta(h_i)
\la(g_i) = \beta(h_i). \]

\noindent The proof can now be completed:
\begin{eqnarray*}
\sg(\bfh) & = & \sum_{\beta \in \Pi',\ \la \in \Pi_p} \prod_{i \in J}
\beta(g_i) \la(h_i) \cdot \prod_{i \notin J} \beta(h_i)\\
& = & \sum_{\beta \in \Pi',\ \la \in \Pi_p} \beta(g_J h_{\notj}) \la(h_J)
= \Sigma_{\Pi'}(g_J h_{\notj}) \cdot \sgp(h_J)
\end{eqnarray*}

\noindent as claimed.
\end{proof}

We compute $\sgp(\bfh)$ in a later section. For now, we mention how to
compute the other factor.

\begin{theorem}\label{Tabel}
Suppose $\ch(R) \nmid |\Pi|$, and no $h_i$ is pseudo-primitive with
respect to (the abelian group of weights) $\Pi$. For each $i$, let $f_i
\in k(R)$ denote $\beta(h_i) / (\beta(g_i) - \beta(g'_i))$ for some
$\beta \in N_{h_i} \cap \Pi$. Also define
$S = \{ I \subset [n] : g_I g'_{\noti} \in \ggh \}$.
Then
$\displaystyle \sg(\bfh) = (-1)^n |\Pi| \prod_{i=1}^n f_i \cdot \sum_{I
\in S} (-1)^{|I|}$.
\end{theorem}

\noindent Note that if $n=1$, this is a special case of the first two
parts of Theorem \ref{T3} above. Moreover, if some $h_i$ {\it is}
pseudo-primitive with respect to $\Pi$, then $\sg(\bfh) = 0$ from
Proposition \ref{P7} above.

\begin{proof}
Let us fix some generators $\beta_1, \dots, \beta_k$ of $\Pi$ (by the
structure theory of finite abelian groups), so that $\Pi =
\bigoplus_{j=1}^k \Z \beta_j$. Then (by assumption,) for each $i$ there
is at least one $j$ so that $\beta_j(g_i) \neq \beta_j(g_i)^{-1}$.

Now fix $i$, and compute $\beta(h_i)$ for arbitrary $\beta \in \Pi$.
Suppose $N'_i$ indexes the set of $\beta_j$'s that are in $N_{h_i}$; then
write $\beta = \beta' + \sum_{j \in N'_i} r_j \beta_j$ for some $r_j \geq
0$ and $\beta' \in \bigoplus_{j \notin N'_i} \Z \beta_j$.

Note by Proposition \ref{P7} for $\Pi \leftrightarrow \bigoplus_{j \notin
N'_i} \Z \beta_j$ that $\beta'(h_i) = 0$.  So if $\beta'' = (\beta')^{-1}
\beta \in \Pi$, then
\[ \beta(h_i) = \beta'(g_i) \beta''(h_i) + \beta'(h_i) \beta''(g_i) =
\beta'(g_i) \beta''(h_i). \]

\noindent It remains to compute the last factor above. This is done using
the following claim.\medskip

\noindent {\bf Claim.} Say $\beta'' = \sum_{j \in N'_i} r_j \beta_j$.
Then $\displaystyle \beta''(h_i) = (\beta''(g_i) - \beta''(g'_i)) f_i$,
where $f_i$ is defined in the statement of the theorem.\medskip

\noindent {\it Proof of the theorem, modulo the claim.} By the claim,
\[ \beta(h_i) = \beta'(g_i)(\beta''(g_i) - \beta''(g'_i))f_i =
(\beta(g_i) - \beta(g'_i))f_i \]

\noindent (for all $i$) since $\beta'(g_i) = \beta'(g'_i)$ by
pseudo-primitivity. Using the notation that $\bbr = \sum_{j=1}^k r_j
\beta_j \in \Pi$, one can compute $\sg(\bfh)$ to be
\begin{align*}
= &\ \sum_{\bf r} \prod_{i=1}^n (\bbr(g_i) - \bbr(g'_i))f_i =
\prod_{i=1}^n f_i \cdot \sum_{\bf r} \sum_{I \subset [n]} (-1)^{n-|I|}
\bbr \left( \prod_{i \in I} g_i \prod_{j \notin I} g'_j \right)\\
= &\ (-1)^n \prod_{i=1}^n f_i \cdot \sum_{I \subset [n]} (-1)^{|I|} \sg
(g_I g'_{\noti})
= (-1)^n \prod_{i=1}^n f_i \cdot \sum_{I \subset [n]} (-1)^{|I|}
\delta_{I \in S} |\Pi|
\end{align*}

\noindent by Proposition \ref{P5}, where the last $\delta$ is 1 if $I \in
S$, and 0 otherwise.
\end{proof}\medskip

The proof is completed by showing the claim.

\begin{proof}[Proof of the claim]
By equation \eqref{E2}, and Lemma \ref{Llevel},
\[ \beta_j^{* r_j}(h_i) = (\beta_j(g_i)^{r_j} - \beta_j(g'_i)^{r_j}) f_i,
\qquad \forall j \in N'_i. \]

\noindent Suppose without loss of generality, that we relabel the set $\{
\beta_j : j \in N'_i \}$ as $\{ \beta_1, \dots, \beta_m \}$ (i.e.,
relabel the generators $\beta_j$ of $\Pi$ so that these are before the
others). Now compute the expression using the above equation:
\[ \beta''(h_i) = \left( \sum_{j=1}^m r_j \beta_j \right) (h_i)
= \sum_{j=1}^m \prod_{l<j} \beta_l(g_i)^{r_l} \cdot (\beta_j(g_i)^{r_j} -
\beta_j(g'_i)^{r_j}) f_i \cdot \prod_{l>j} \beta_l(g'_i)^{r_l}, \]

\noindent and this telescopes to
$f_i \cdot \prod_j \beta_j(g_i)^{r_j} - f_i \cdot \prod_j
\beta_j(g'_i)^{r_j} = (\beta''(g_i) - \beta''(g'_i))f_i$, as claimed.
\end{proof}

\section{Products of pseudo-primitive elements - positive
characteristic}\label{Spseudo}

We now mention results for pseudo-primitive elements $h_i$ (and not
necessarily abelian $\Pi$) in prime characteristic; note that for almost
all characteristics (including zero), Proposition \ref{P7} above says
that $\sg(\bfh) = 0$. Before considering the positive case, we need a
small result.

\begin{lemma}\label{LBern}
Given $f \in \N$ and a prime $p>0$, define $\varphi_p(f) = \varphi_p(f)
:= \sum_{i=0}^{p-1} i^f$. If $f>0$, then $\varphi(f) \not\equiv 0 \mod p$
if and only if $(p-1)|f$, and in this case, $\varphi(p-1) = p-1 \equiv -1
\mod p$.
\end{lemma}

\begin{proof}
Let $g$ be any cyclic generator of (the finite cyclic group) $(\Z /
p\Z)^\times$. Then $\sum_{i=1}^{p-1} i^f \equiv \sum_{j=1}^{p-1} g^{jf}
\mod p$. Now if $(p-1) | f$, then each summand is 1, which yields $p-1
\mod p$. Otherwise, $g^f$ is not 1, and its powers add up to 0 (by the
geometric series formula).
\end{proof}

\subsection{Preliminaries}

Recall that a pseudo-primitive element is any $h \in H_{g,g'}$ so that
$g^{-1}g' \in \ggh$.
%
%
Some terminology is now needed. Note by Hall's theorems, that a finite
group $\Phi$ is solvable if and only if it contains Hall subgroups of all
possible orders (e.g., see \cite[\S 11]{AB}).  So if $|\Phi| = p^k \cdot
m$ with $p \nmid m$, let $\Phi_m$ be any Hall subgroup of order $m$.

\begin{defn}
($p > 0$ a fixed prime.) Given a finite solvable group $\Phi$, denote by
$\Phi_p, \Phi'$ respectively, any $p$-Sylow subgroup and any Hall
subgroup of order $|\Phi| / |\Phi_p|$. (From above, we mean $\Phi' =
\Phi_m$.)
\end{defn}

\noindent For the rest of this section, $\ch(R) = p>0$; also fix $n$, the
number of $h_i$'s. (Recall Assumption \ref{D51}.)

\begin{prop}\label{P4}
Suppose, given skew-primitive $h_i \in H_{g_i, g'_i}$ for $1 \leq i \leq
n$, that $\Phi \subset \Gamma / \Gamma'$ is a finite solvable subquotient
of $\Gamma$ (as in a previous section) with respect to which every $h_i$
is pseudo-primitive. Then $\Sigma_{\Phi}(\bfh) = \sgphi(\bfg) \cdot
\sgphip(\bfh)$.
\end{prop}

\begin{proof}
The first claim is that the ``set-product" $\Phi' \Phi_p := \{ \beta *
\la : \beta \in \Phi', \la \in \Phi_p \}$ equals the entire group $\Phi$.
\comment{
For if not, then some two elements $\beta * \la, \beta' * \la'$ are
equal. But then $\beta^{-1} \beta' = \la (\la')^{-1}$, whence this common
element has order dividing both $p^k$ and $m$. Thus $\beta = \beta', \la
= \la'$, so the multiplication map is injective, and between two finite
sets of the same order.
}
Next, if $\beta \in \Phi'$, then $\beta(h) = 0$ for any pseudo-primitive
$h \in H_{g,g'}$, because if $|\Phi'| = m \not\equiv 0 \mod p$, then
$\beta^{*m} = \vi \in \Gamma' \subset \Gamma_h$, whence $0 = \vi(h) =
\beta^{*m}(h) = m \beta(g)^{m-1} \beta(h)$. Therefore,
\begin{eqnarray*}
\Sigma_{\Phi}(\bfh) & = & \sum_{\beta \in \Phi',\ \la \in \Phi_p}
\prod_{i=1}^n (\beta * \la)(h_i) = \sum_{\beta \in \Phi',\ \la \in
\Phi_p} \prod_{i=1}^n \beta(g_i)\la(h_i)\\
& = & \sum_{\beta \in \Phi'} \prod_{i=1}^n \beta(g_i) \cdot \sum_{\la \in
\Phi_p} \prod_{i=1}^n \la(h_i) = \sgphi(\bfg) \sgphip(\bfh)
\end{eqnarray*}

\noindent as claimed.
\end{proof}

\begin{remark}
The above proposition thus holds for {\it any} group $\Phi$, such that
some $p$-Sylow subgroup $\Phi_p$ has a complete set of coset
representatives, none of whom has order divisible by $p$. Obvious
examples are abelian groups or groups of order $p^a q^b$ for primes $p
\neq q$ (but these are solvable by Burnside's Theorem).
\comment{
Moreover, J.~Alperin mentioned the following result in \cite{Alp}:
\begin{prop}
Suppose a finite group $G$ has order $\prod_i p_i^{a_i}$ for distinct
primes $p_i$. Then $G$ is solvable if and only if there exist $p_i$-Sylow
subgroups for each $i$, whose ``set-products" are all Hall subgroups of
$G$.
\end{prop}
}

Also note that the above sum is independent of the choices of $\Phi_p,
\Phi'$.
\end{remark}\medskip

\noindent The next result is crucial in computing $\sg(\bfh)$, and uses
subquotients of $\Pi$.

\begin{theorem}\label{T5}
Given a finite subgroup of weights $\Pi \subset \Gamma = \Gamma_H$,
suppose $h_i \in H_{g_i, g'_i}$ is pseudo-primitive with respect to $\Pi$
for all $1 \leq i \leq n$. Define $\Phi = \Pi_{ab} := \Pi / [\Pi, \Pi]$.
Then,
\begin{equation}\label{Ephi}
\sg(\bfh) = |[\Pi, \Pi]| \cdot \sgphi(\bfg) \cdot \sgphip(\bfh).
\end{equation}
\end{theorem}

\noindent For instance, if every $h_i$ was almost primitive, then
$\sgphi(\bfg) = [\Phi : \Phi_p]$.

\begin{proof}
At the outset, note that $\sgphi(\bfg)$ and $\sgphip(\bfh)$ make sense
because of Corollary \ref{C1} and Proposition \ref{Pstab} above. Now, the
proof is in two steps; each step uses a previously unused result
above.\medskip

\noindent {\bf Step 1.} We claim that $\sg(\bfh) = |[\Pi, \Pi]|
\Sigma_{\Phi}(\bfh)$. This follows immediately from Lemma \ref{L1}, where
$h,\Gamma'',\Pi^\circ$ are replaced by $\bfh, [\Pi, \Pi], \Pi$
respectively.

The only thing to check is that the above replacements are indeed valid.
Since $[\Pi, \Pi]$ is normal in $\Pi$, it suffices to check that $[\Pi,
\Pi] \subset \Gamma_{\bfh}$. But this follows from Corollary \ref{C1} and
Proposition \ref{Pstab} above.\medskip

\noindent {\bf Step 2.} The proof is now complete by invoking Proposition
\ref{P4} above.
\end{proof}

We conclude the preliminaries with one last result - for {\it
skew}-primitive elements in general.

\begin{prop}
If $\Pi_p \ncong (\Z / p\Z)^k$ for any $k>0$, then $\Sigma_{\Pi_p}(\bfh)
= 0$.
\end{prop}

\comment{
\noindent Therefore if we want to compute $\Sigma_{\Pi_p}(\bfh)$ for
skew-primitive $h_i$'s, the only interesting case is when $\Pi_p \cong
(\Z / p\Z)^k$ for some $k$ (by Theorem \ref{Tskew}); we address this case
below.
}

\begin{proof}
By Theorem \ref{Tskew}, the $h_i$'s are almost primitive with respect to
$\Pi_p$. Now invoke equation \eqref{Ephi} above, replacing $\Pi$ by
$\Pi_p$. Now, if $\Pi_p$ is not abelian, then $|[\Pi_p, \Pi_p]| > 1$,
hence is a power of $p$, whence the right-hand side vanishes. Next, if
$\Pi_p$ is abelian, but contains an element of order $p^2$, then
$\sgp(\bfh) = 0$ by Theorem \ref{Tskew} again. Therefore $\Pi_p \cong (\Z
/ p\Z)^k$ for some $k$. If $k=0$, then $\Pi_p = \{ \vi \}$, and
$\vi(\bfh) = 0$.
\end{proof}

\subsection{The main results - pseudo-primitive elements}

The following result now computes $\sg(\bfh)$ (for pseudo-primitive
$h_i$'s) in most cases in prime characteristic that are ``nonabelian".
For the ``abelian" case, we appeal to Theorem \ref{T8} below - and
mention at the outset, that it is only for {\it almost} primitive (and
not merely pseudo-primitive) elements, that we obtain a much clearer
picture - as its last part shows.

\begin{theorem}\label{T6}
Suppose $\ch(R) = p \in \N$, $\Pi_p$ is any (fixed) $p$-Sylow subgroup of
$\Pi$, and every $h_i$ is pseudo-primitive (with respect to $\Pi$).
\begin{enumerate}
\item $\sg(\bfh) = 0$ if $\Pi_p$
\begin{enumerate}
\item is trivial,
\item contains an element of order $p^2$, or
\item intersects $[\Pi, \Pi]$ nontrivially.\medskip

This last part includes the cases when $\Pi_p$
\item is not abelian,
\item does not map isomorphically onto (some) $\Phi_p$, via (the
restriction of) the quotient map $\pi : \Pi \twoheadrightarrow \Phi = \Pi
/ [\Pi, \Pi]$, or
\item has size strictly greater than $\Phi_p$.
\end{enumerate}

\item Otherwise $\Pi_p \cong \Phi_p \cong (\Z / p\Z)^k$ for some $k>0$,
and then $\sg(\bfh) = |[\Pi, \Pi]| \cdot \sgphi(\bfg) \cdot \sgp(\bfh)$.
\end{enumerate}
\end{theorem}

\begin{remark}
Any finite abelian group of exponent $p$ is of the form $(\Z / p\Z)^k$,
hence one part of the second statement is clear. Moreover, every
subquotient of such a group is of the same form. Finally (especially when
all of the $h_i$'s are almost primitive with respect to $\Pi$), the cases
that remain reduce to computing $\sgp(\bfh)$, and when $\Pi_p \cong (\Z /
p\Z)^k$; this is addressed below.
\end{remark}

\renewcommand{\theenumi}{\alph{enumi}}
\begin{proof}
The second part follows from the first part, the remarks above, and
Equation \eqref{Ephi}. We now show the first part.
\begin{enumerate}
\item[(a)] If $\Pi_p$ is trivial, then $p \nmid |\Pi|$, and we are done
by Proposition \ref{P7}.

\item[(b)] This has been done in Theorem \ref{Tskew} above.

\item[(c)] Now suppose that $[\Pi, \Pi] \cap \Pi_p \neq \emptyset$. Then
$[\Pi, \Pi]$ contains an element of order $p$, whence $p$ divides $|[\Pi,
\Pi]|$. Now use equation \eqref{Ephi}.\medskip

\noindent It remains to show how this last includes the remaining
cases.

\item[(d)] First, if $\Pi_p$ is nonabelian, then $[\Pi_p, \Pi_p]$ is a
nontrivial subgroup of the $p$-group $\Pi_p$. In particular, $\Pi_p$
intersects $[\Pi, \Pi]$.\medskip

\item[(e)] Next, note that $\pi(\Pi_p)$ is a $p$-group in $\Phi$, and
$|\Pi_p| \geq |\Phi_p|$ (since $|\Phi|$ divides $|\Pi|$). Hence $\Pi_p$
does not map isomorphically onto (some) $\Phi_p$ if and only if $\pi$ is
not one-to-one on $\Pi_p$. But then $[\Pi, \Pi]$ intersects $\Pi_p$.

\item[(f)] Finally, if $|\Pi_p| > |\Phi_p|$, then $\Pi_p$ cannot map
isomorphically onto $\Phi_p$, so we are done by the preceding paragraph.
\end{enumerate}
\end{proof}

We conclude by analyzing $\Sigma_{\Pi_p}(\bfh)$. Note that the results
below that pertain only to $\Sigma_{\Pi_p}(\bfh)$ are applicable in
general to all skew-primitive $h_i$'s, by Theorem \ref{Tskew} above.

\renewcommand{\theenumi}{\arabic{enumi}}
\begin{theorem}\label{T8}
Suppose $\ch(R) = p \in \N$, $\Pi_p$ is any (fixed) $p$-Sylow subgroup of
$\Pi$, and every $h_i \in H_{g_i, g'_i}$ is pseudo-primitive (with
respect to $\Pi$) for all $1 \leq i \leq n$. Suppose moreover that $\Pi_p
\cong (\Z / p\Z)^k$.
\begin{enumerate}
\item $\g(g_i) = 1$ for all $i$ and $\g \in \Pi_p$. In particular, $g_i,
\bfg \in G_{\Pi_p}(H)$.

\item If $k > n$, then $\sg(\bfh) = \sgp(\bfh) = 0$.

\item If $k=n$, then
$\displaystyle \sgp(\bfh) = \binom{p}{2}^k \cdot \perm(A)$,
where $A$ is the matrix given by $a_{ij} = \g_j(h_i)$, the
$\g_j$'s form a $\Z / p\Z$-basis of $\Pi_p$, and $\perm$ is the {\em
matrix permanent}:
\[ \perm(A_{n \times n}) = \sum_{\sigma \in S_n} \prod_{i=1}^n a_{i,
\sigma(i)}. \]

\noindent In particular, $\sgp(\bfh) = 0$ unless $p=2$, in which case
$\sg(\bfh) = \sgp(\bfh) = \det A$.
\comment{
Moreover, the sum actually looks like a determinant, but without the
signs (but we can add the appropriate signs modulo 2); {\it this}
expression is called the {\it permanent} of the matrix $A$.
}

\item If $\sgp(\bfh) \neq 0$, then $(p-1) | n$ and $0< k \leq n / (p-1)$,
and then $\sgp(\bfh)$ can take any value $r \in R$. (If $k=n$ and $p=2$,
then $r \neq 0$.)\medskip
\end{enumerate}
\end{theorem}

\begin{remark}\label{R711}
This result is independent of the chosen $p$-Sylow subgroup
$\Pi_p$, as well as the choices of generators $\g_j$. It generalizes
Proposition \ref{Punproved} above, in the special case $p=2$.
\comment{
\begin{enumerate}
\item Let us look at the special case $p=2$ (and also focus on $n=1,\
g,g' \in \ggh$ as considered in Proposition \ref{Punproved} above). If
$\Pi_2$ denotes any fixed 2-Sylow subgroup, and has an element of order
4, then $\sg(\bfh) = 0$. Otherwise every nontrivial element in $\Pi_2$
has order 2, whence
$\Pi_2$ is a finite abelian group with exponent 2, i.e., $\Pi_2 \cong (\Z
/ 2\Z)^k$ for some $k$. Theorem \ref{T8} then tells us that $\sg(\bfh) =
0$ if $k > n$. (If $n=1$, then this happens when 4 divides $|\Pi_2|$.)

Moreover, since $p=2$, the last part says (if $\bfg \in \ggh$ as well)
that $\sg(\bfh) = \sgp(\bfh)$ can take any value for any $k \leq n =
n/(p-1)$. For instance, if $n=1$, then $\sgp(h) = \g(h)$ (the determinant
- or permanent - of a $1 \times 1$ matrix) for any $\g$ of order 2 (i.e.,
a basis of some 2-Sylow subgroup).
\end{enumerate}
}
\end{remark}\medskip

\noindent The rest of this section is devoted to proving the above
result. First, suppose there exists a subgroup $\Pi'_p \cong (\Z /
p\Z)^k$ of $\Pi$ (so $\Pi'_p \subset \Pi_p$ in general). Choose a set of
coset representatives $\B$ for $\Pi'_p$ in $\Pi$, and write
\[ \sg(\bfh) = \sum_{\beta \in \B,\ \g \in \Pi'_p} (\g * \beta)(\bfh) =
\sum_{\beta \in \B,\ \g \in \Pi'_p} \prod_{i=1}^n (\g * \beta)(h_i). \]

Recall that every element of $\Pi'_p$ is $\g_{a, I} := \sum_{j \in I} a_j
\g_j$ (with $\g_i$ as above), for some subset $I$ of $[k] := \{ 1, 2,
\dots, k \}$, and some $|I|$-tuple $a = (a_j)_{j \in I}$ of elements of
$(\Z / p\Z)^\times$. Recall, moreover, that we had previously defined
$g_I, g'_I, h_I$ for $I \subset [n]$.\medskip

\noindent We also need the following lemma, that is proved by
computations using Proposition \ref{P8}.

\begin{lemma}
If $h \in H_{g,g'}$ is pseudo-primitive with respect to $\Pi$, and given
$\beta \in \B, \g_{a,I} \in \Pi_p$ as above, one has
\[ (\beta * \g_{a,I})(h) = \prod_{j=1}^k \g_j(g)^{a_j} \cdot \left[
\beta(h) + \beta(g) \sum_{j=1}^k a_j \g_j(g^{-1}h) \right]. \]
\end{lemma}

\comment{
\begin{proof}
We sketch this (purely computational) proof: by Proposition \ref{P8}, we
know $\g_j^{*a_j}(h)$ for all $j$, and since $h$ is pseudo-primitive, we
replace all $g'$'s by $g$'s. Therefore denoting $\g_{a',I'} := \sum_{j<k}
l_j \g_j$, we inductively get that
\begin{eqnarray*}
\g_{a,I}(h) = (\g_{a',I'} * l_k \g_k)(h) & = & \g_k(g)^{l_k} \cdot
\prod_{j<k} \g_j(g)^{l_j} \sum_{j<k} l_j \g_j(g^{-1}h)\\
& + & \prod_{j<k} \g_j(g)^{l_j} \cdot l_k \g_k(g)^{l_k - 1} \g_k(h)
\end{eqnarray*}

\noindent which shows the result for $\beta = \vi$. Now we just convolve
(convolute?) with $\beta$ to get the given result.
\end{proof}\medskip
}

The key observation now, is that the only ``monomials" that occur in the
product $\prod_{i=1}^n (\beta * \g_{a,I})(h_i)$ are of the form
$\beta(\bfg g_{I_0}^{-1} h_{I_0}) \cdot \prod_{j \in I} \g_j(h_{I_j})$,
where $\coprod_j I_j \coprod I_0 = [n]$, and $I_j \subset I$
for all $j$. The coefficient of such a monomial in this particular
summand, is $\prod_{j \in I} a_j^{|I_j|} \g_j(\bfg)^{a_j}$ by the lemma
above. Moreover, every such monomial occurs at most once inside each
$(\beta * \g_{a,I})(\bfh)$.\medskip

\noindent The crucial fact that proves Theorem \ref{T8} above, is the
following\medskip

\noindent {\bf Key claim.} The coefficient of $\beta(\bfg g_{I_0}^{-1}
h_{I_0}) \prod_{j \in I} \g_j(h_{I_j})$ in $\sum_{\g \in \Pi'_p} (\beta *
\g)(\bfh)$, equals\hfill\break
$p^{k-|I|} \prod_{j \in I} \varphi(|I_j|)$ (where $\varphi(0) := 0$).

\begin{proof}[Proof of the key claim]
(Note that $\varphi(f)$ was defined in Lemma \ref{LBern}.)
A monomial of the desired form occurs in precisely those
$(a',I')$-summands, so that $I' \supset I$.
Moreover, all such summands can be split up into a disjoint union over
all $a \in ((\Z / p\Z)^\times)^{|I|}$, with each disjoint piece
containing all $(a', I')$ so that $I' \supset I$ and the $I$-component of
$a'$ is $a$.

Such a piece contains exactly $p^{k-|I|}$ elements (and hence exactly
that number of copies of the monomial with this selfsame coefficient).
Each of these ``extra" $[k] \setminus I$ factors contributes a
$\beta(h_i)$, which gives $\beta(h_{I_0})$.

Moreover, there is one contribution for each $a \in ((\Z /
p\Z)^\times)^{|I|}$, and it is $\prod_{j \in I} a_j^{|I_j|}
\g_j(\bfg)^{a_j} \cdot \beta(\bfg g_{I_0}^{-1})$, since the argument for
the $\beta$-factor here is precisely $\prod_{j=1}^k g_{I_j}$. Moreover,
$\g_j(\bfg) = 1$ by Theorem \ref{Tskew} above.

Summing over all possible tuples $a \in ((\Z / p\Z)^\times)^{|I|}$, the
coefficient (apart from the $\beta$-part) is
\[ p^{k-|I|} \sum_a \prod_{j \in I} a_j^{|I_j|} = p^{k-|I|} \prod_{j \in
I} \sum_{a_j = 1}^{p-1} a_j^{|I_j|}, \]

\noindent and this equals $p^{k-|I|} \prod_{j \in I} \varphi(|I_j|)$ as
desired, because the only problem may occur when some $|I_j| = 0$. But
then $|I| < k$, so
\[ p^{k-|I|} \sum_{a_j=1}^{p-1} a_j^0 = 0 \sum_{a_j=1}^{p-1} a_j^0 = 0
\sum_{a_j=1}^p a_j^0 = p^{k-|I|} \varphi(0). \]
\end{proof}

\begin{proof}[Proof of Theorem \ref{T8}]\hfill
\begin{enumerate}
\item This is from Theorem \ref{Tskew} above.
\end{enumerate}\medskip

\noindent Now set $\Pi_p = \Pi'_p$.
We first note from the key claim that if $I_0$ is nonempty, or any $I_j$
is empty, then the coefficient of that particular monomial vanishes -
because $\ch(R) = p$ and $\varphi(0) = 0$.

\begin{enumerate}
\setcounter{enumi}{1}
\item Suppose $k>n$. Then at least one $I_j$ must be empty in every
monomial above, by the Pigeonhole Principle, and we are done.\medskip

\item If $k=n$, then the only monomials that have a nonzero contribution
to the sum $\Sigma_{\Pi'_p}(\bfh)$ must correspond to empty $I_0$ and
singleton $I_j$'s (since $\coprod_{j=1}^{k} I_j = [n]=[k]$). In other
words,
\comment{
\item Now choose $k=n$ in $\Pi'_p$ (so $|\Pi_p| \geq p^n$). Then from the
key claim, first observe that the only monomials that have nonzero
coefficient satisfy $|I| = k = n$ here. Thus, $\beta$ has no contribution
(i.e., $|I_0| = 0$), whence the sum $\sum_{\g \in \Pi'_p} (\beta *
\g)(\bfh)$ is independent of $\beta$. Hence $\sg(\bfh) = [\Pi : \Pi'_p]
\Sigma_{\Pi'_p}(\bfh)$. This analysis works if we replace $\Pi$ by
$\Pi_p$ as well, so $\sg(\bfh) = \sgp(\bfh) = 0$ (modulo $p$) if $|\Pi_p|
> |\Pi'_p|$.

Next, if $|\Pi_p| = |\Pi'_p| = p^n$, the monomials with possibly nonzero
terms must be of the form $\prod_{j=1}^n \g_j(h_{I_j})$, with $\coprod_j
I_j = [n]$. This implies that each $I_j$ is a singleton, giving us
}
$\sigma \in S_n : j \mapsto i_j\ \forall j$. Moreover, the coefficient of
such a monomial is $p^0 \prod_{j=1}^n \varphi(1)$, and these monomials
all add up to give the matrix permanent, as claimed. The rest of the
statements are now easy to see.\medskip

\item In this part, we are only concerned with $\sgp(\bfh)$, so that
$\beta$ does not contribute here either (so $I_0 = \emptyset$ and $[n] =
\coprod_{j \in I} I_j$).

From the key claim and Lemma \ref{LBern} above, observe that if some
monomial has a nonzero contribution, then $(p-1)$ divides $|I_j|$ for all
$j$, and $I = [k]$. In particular, $(p-1)$ divides $\sum_{j \in I} |I_j|
= n$, and
\[ n = \sum_{j \in I} |I_j| = \sum_{j=1}^k |I_j| \geq \sum_{j=1}^k (p-1)
= k(p-1), \]

\noindent whence $k \leq n/(p-1)$. Moreover, $\sgp(\bfh) = \vi(\bfh) = 0$
if $k=0$.

It remains to present, for each $0 < k \leq n/(p-1)$ and (nonzero) $r \in
R$, an example of $(H, \Pi = \Pi_p)$, so that $\sg(\bfh) = \sgp(\bfh) =
r$. This example is analyzed in the next section.
\end{enumerate}
\end{proof}

\section{Example: Lie algebras}

Suppose $H = U(\mfg)$ for some Lie algebra $\mfg$ (say over $\comp$).
Then any weight $\mu \in \Gamma_H$ kills $[\mfg,\mfg]$, hence belongs to
$(\mfg/[\mfg,\mfg])^*$. Let us denote $\mfgone := \mfg / [\mfg,\mfg]$.
Conversely, any element of the set above, is a weight of $H$, using
multiplicativity and evaluating it at the projection down to the quotient
$\mfgone$.
Thus, $\Gamma_H$ is the dual space (under addition) of the abelianization
$\mfgone$ of $\mfg$. Hence we now examine what happens in the case of an
($R$-free) abelian Lie algebra $\mfh$.\medskip

In this case, we have the free $R$-module $\mfh = \oplus_i R h_i$ with
the trivial Lie bracket, and $H = U(\mfh) = \Sym(\mfh)$. Thus, $H$
inherits the usual Hopf algebra structure now (i.e., $\dd(h_i) = 1
\otimes h_i + h_i \otimes 1,\ S(h_i) = -h_i,\ \vi(h_i) = 0\ \forall i$).

First, $(\Gamma_H,*) = (\mfh^*, +)$. By Proposition \ref{P7}, if $\ch(R)
\nmid |\Pi|$, then $\sg(\bfh) = 0$ for all products $\bfh$ of primitive
elements in $H$ (and hence for all $\bfh$ in the augmentation ideal $\mfh
U(\mfh)$ of $H$). Thus, the only case left to consider is when $\ch(R) =
p > 0$. But then $(\mfh^*, +)$ is a $\Z / p\Z$-vector space, so every
finite subgroup $\Pi$ is of the form $\Pi = \Pi_p \cong (\Z / p\Z)^k$ for
some $k$. Moreover, Theorem \ref{T3} and (the last part of) Theorem
\ref{T8} provide more results in this case.\medskip

We therefore conclude the example (and the proof of the theorem above) by
analyzing the computation of $\sg(\bfh)$ for $\bfh = h_1 \dots h_n$. For
any (nonzero) $r \in R$, we produce such a finite subgroup $\Pi = \Pi_p
\cong (\Z / p\Z)^k$, so that $0 < k \leq n/(p-1)$ and $\sg(\bfh) =
r$.\medskip

\noindent {\it Construction:} Given $k$, partition $[n]$ into $k$
disjoint nonempty subsets $[n] = \coprod_{j=1}^k I_j$, reordered so that
$I_1 = \{ 1, \dots, n- (k-1)(p-1) \}$, and so that $|I_j| = p-1$ for all
$j>1$.
For each $1 \leq j \leq n$, define $\g_j \in \mfh^* = \Gamma_H$ as
follows: $\g_1(h_1) = r,\ \g_j(h_i) = 1$ if $i \neq 1 \in I_j$, and
$\g_j(h_i) = 0$ otherwise. (One verifies that the $\g_i$'s thus defined
are indeed linearly independent over $k(R)$, hence over $\Z / p\Z$ as
well, but for this, one needs that $r \neq 0$ if $n=k, p=2$.) Thus for
any $K \subset [k],\ \g_K(h_i) := \sum_{j \in K} \g_j(h_i)$ vanishes
unless $i \in \cup_{j \in K} I_j$.

Now evaluate $\sg(\bfh) =$ $\sum_{K \subset [k], a} \prod_{i=1}^n
\g_{a,K}(h_i)$, where
$\Pi := \sum_{i=1}^n \Z \g_i = \bigoplus_{i=1}^n (\Z / p\Z) \g_i$.
By the key claim in the previous section, the only monomials $\prod_{j
\in I} \g_j(h_{I'_j})$ that do not vanish are for $|I| = k$, and with
$(p-1)$ divides $|I'_j|$ for all $j$.
Moreover, $\g_j(h_i)$ is zero except when $i \in I_j$, so there is only
one type of monomial remaining: $\prod_{j \in I} \g_j(h_{I_j})$. (Note
that this satisfies the earlier condition: $(p-1)$ divides $|I_j|$ for
all $j$.)

Moreover, by the key claim in the preceding section, the coefficient of
this monomial, which itself equals $r \cdot \prod_{i=2}^n 1 = r$, is
$\prod_{j=1}^k \varphi(|I_j|)$, and by Fermat's Little Theorem,
$\varphi(|I_j|) = p-1 = -1\ \forall j$ (in characteristic $p$). We
conclude that $\sg(\bfh) = \sgp(\bfh) = (-1)^k r$, whence we are done
(start with $r' = (-1)^k r$ to get $r$).

\section{Example: Degenerate affine Hecke algebras of reductive type with
trivial parameter}

In this section, we apply the general theory above, to a special case,
wherein a finite group acts on a vector space (or free $R$-module in our
case), with the group and the module corresponding to the Weyl group and
the Cartan subalgebra (actually, its dual space) respectively, of a
reductive Lie algebra. We use the $\Z$-basis of simple roots (and any
$\Z$-basis for the center), to try and compute the value of $\sg(\bfh)$.

\subsection{Hopf algebras acting on vector spaces}

We will consider special cases of the following class of Hopf algebras.
Suppose that a cocommutative $R$-Hopf algebra $H$ acts on a free
$R$-module $V$; denote the action by $h \cdot v$ for $h \in H, v \in V$.
Then $H$ also acts on $V^*$ by: $\tangle{h \cdot \la, v} := \tangle{\la,
S(h) \cdot v}$.

Now consider the $R$-algebra $A$ generated by the sets $H$ and $V$, with
obvious relations in $H$, and the extra relations $vv' = v'v,\ \sum
\one{h} v S(\two{h}) =: \ad h(v) = h \cdot v$ for all $h \in H$ and $v,v'
\in V$. Note that the relation $\ad h(v) = h \cdot v$ can be rephrased,
as the following lemma shows.

\begin{lemma}
Suppose some $R$-Hopf algebra $H$ acts on a free $R$-module $V$, and an
$R$-algebra $B$ contains $H,V$. Then the following relations are
equivalent (in $B$) for all $v \in V$:
\begin{enumerate}
\item $\sum \one{h} v S(\two{h}) = h \cdot v$ for all $h \in H$.
\item $h v = \sum (\one{h} \cdot v) \two{h}$ for all $h \in H$.
\end{enumerate}

\noindent If $H$ is cocommutative, then both of these are also equivalent
to:

\begin{enumerate}
\setcounter{enumi}{2}
\item $vh = \sum \one{h} (S(\two{h}) \cdot v)$ for all $h \in H$.
\end{enumerate}

\noindent Moreover, if this holds, then any unital subalgebra $M$ of $B$
that is also an $H$-submodule (via $\ad$), is an $H$-(Hopf-)module
algebra under the action
\[ h \cdot m := \ad h(m) = \sum \one{h} m S(\two{h})\ \forall h \in H,\ m
\in M. \]
\end{lemma}

\noindent (The proof is straightforward.) For instance, one can take $M =
B$ or $H$ - or in the above example of $A$, consider $M = \Sym_R
V$.\medskip

It is straightforward (but perhaps tedious) to check that $A$ is a Hopf
algebra with the usual operations: on $H$, they restrict to the Hopf
algebra structure of $H$, and $V$ consists of primitive elements.

By the above lemma, if $H$ is $R$-free, then the ring $A$ is an $R$-free
$R$-Hopf algebra, with $R$-basis given by $\{ h \cdot m \}$, where $h \in
H$ and $m$ run respectively over some $R$-basis of $H$, and all
(monomial) words (including the empty word) with alphabet given by an
$R$-basis of $V$.
It has the subalgebras $H$ and $\Sym_R(V)$, and is called the {\it smash
product} $H \ltimes \Sym_R V$ of $H$ and $\Sym_R V$.\medskip

We now determine the weights of $A$. Denote by $\Gamma_H$ the group of
weights of $H$ (under convolution). One can now use Proposition \ref{Pq}
to prove:

\begin{prop}\label{Pweights}
The weights $\Gamma_A$ of $A$ form a group, which is the Cartesian
product $\Gamma_H \times V^*_\vi$, with convolution given by
\[ (\nu_1, \la_1) * (\nu_2, \la_2) = (\nu_1 \nu_2, \la_1 + \la_2) =
(\nu_1 *_H \nu_2, \la_1 *_V \la_2) \]

\noindent for $\nu_i \in \Gamma_H,\ \la_i \in V^*_\vi$. (Here, $V^*_\vi$
is the $\vi$-weight space of the $H$-module $V^*$.)
\end{prop}

\comment{
\begin{proof}
Given $\g \in \Gamma$, its restrictions to $A$ and $\Sym_R(V)$ are also
algebra maps, or weights, of the respective subalgebras. From above, the
sets of weights of these subalgebras are $\Gamma_A$ and $V^*$
respectively. Moreover, these ``restricted" weights extend to the weight
$\g$ of $H$ if and only if $\g(\ad a(v)) = \g(a(v))$, i.e.,
\[ \g(a(v)) = \g(\ad a(v)) = \sum \g(\one{a}) \g(S(\two{a})) \g(v) =
\vi(a) \g(v). \]

\noindent Thus, the condition above for all $a,v$ is equivalent to saying
(note that $A$ is cocommutative, so $S^2 = \id$ on $A$, whence $S$ is a
bijection on $A$) that
\[ a(\g)(v) = \tangle{\g, S(a)(v)} = \mbox{(above) } \vi(S(a)) \g(v) =
\vi^{-1}(a) \g(v) = \vi(a) \g(v) \]

\noindent (since $\vi$ is the identity in $\Gamma$). Thus for all $a$,
$a(\g) - \vi(a) \g$ kills all of $V$.\medskip

Conversely, given $\nu \in \Gamma_A$ and $\la \in V^*_\vi$, one observes
that these can be extended to get a weight $\g = (\nu, \la) \in \Gamma$,
since the defining relation holds:
\begin{eqnarray*}
\g(\ad a(v)) & = & \sum \nu(\one{a}) \nu(S(\two{a})) \la(v) = \vi(a)
\la(v) = (\vi(S(a)) \la)(v)\\
& = & S(a)(\la)(v) = \la(a(v))
\end{eqnarray*}

\noindent and hence the two sets above are equal, as claimed.
\end{proof}
}

\subsection{Degenerate affine Hecke algebras}

Since we work over any commutative unital integral domain $R$, hence we
can generate examples over all $R$, if there exists a lattice in $V$ that
is fixed by $H$, and one considers its $R$-span. Now specialize to the
case when $H = RW$ is the group ring of a Weyl group acting on a Cartan
subalgebra of the corresponding semisimple Lie algebra. Then one uses the
root lattice $Q$ inside $V = \mfh^*$.

We work in slightly greater generality. Given a finite-dimensional
reductive complex Lie algebra $\mfg$, let $W$ be its Weyl group and
$\mfh$ a fixed chosen Cartan subalgebra. Thus $\mfh = \oplus_{i \geq 0}
\mfh_i$, where for $i>0,\ \mfh_i$ corresponds to a simple component
(ideal) of $\mfg$, with corresponding base of simple roots $\dd_i$ and
Weyl group $W_i$, say; and $\mfh_0$ is the central ideal in $\mfg$.

Define $Q_i = \oplus_{\aaa \in \dd_i} \Z \aaa$, the root lattice inside
$\mfh_i^*$, and choose and fix some $\Z$-lattice $Q_0$ inside $\mfh_0^*$.
Now replace $\mfh_i^*$ by $V_i = \mfh_i^* := R \otimes_{\Z} Q_i$, and
$\mfh_i$ by the $R$-dual of $\mfh_i^*$, for all $i \geq 0$. Thus, for the
entire Lie algebra, $\dd = \coprod_{i > 0} \dd_i$ and $W = \times_{i > 0}
W_i$.

Now define $V = \oplus_{i \geq 0} V_i$, whence the previous subsection
applies and one can form the algebra $A = RW \ltimes \Sym_R V$. This is
the {\it degenerate affine Hecke algebra with trivial parameter} (the
parameter is trivial since $wv - w(v)w$ is always zero), of reductive
type. This is a special case of \cite[Definition 1.1]{Ch}, where one sets
$\eta = 0$.\medskip

Before we address the general case, note that there are two types of
$\mfh_i$'s in here: ones corresponding to simple Lie algebras, which we
address first, and the ``central part", which is fixed by $W$ (hence so
is $\mfh_0^*$).

\subsection{The simple case}

The first case to consider is: $V = \mfh^* = R \otimes_{\Z} Q$, for a
{\it simple} Lie algebra. Thus $\dd$ is irreducible, and given $A = RW
\ltimes \Sym_R(\mfh^*)$, $\Gamma_A = \Gamma_W \times \mfh^W$ (because the
condition in Proposition \ref{Pweights} above translates to: $w(\g) =
\vi(w) \g = \g$ for all $w \in W, \g \in \Gamma_A$). Here, $\Gamma_W =
\Gamma_{RW}$.

We now state our main result, using the convention that all roots in the
simply laced cases (types $A,D,E$) are short. The result helps compute
$\sg$ at any element of the $R$-basis $\{ g \cdot m \}$ mentioned in an
earlier subsection.

\begin{theorem}\label{T7}
Suppose $\mfg$ is a complex simple Lie algebra with simple roots
$\Delta$, Weyl group $W$, $V = \mfh^* = R \otimes_{\Z} \Z \Delta$, and $A
= RW \ltimes \Sym_R(\mfh^*)$. As above, let $\Pi \subset \Gamma_A$ be a
finite subgroup of weights. Let $h_1, \dots, h_n \in \mfh^*$.
\begin{enumerate}
\item If $\ch(R) \neq 2$, or $W$ is of type $G_2$, or $W$ has more than
one short simple root, then every weight acts as $\vi = 0$ on $\mfh^*$.
In particular, $\sg = 0$ on $\Sym_R(\mfh^*)$.

\item If $\ch(R) = 2$, then every weight acts as $\vi$ on $W$. Now
suppose also that $W$ is not of type $G_2$, and has only one short simple
root $\aaa_s$, say.

If $\Pi$ has an element of order 4, or $h_i$ has no
``$\aaa_s$-contribution" (i.e., $h_i \in \oplus_{\aaa_s \neq \aaa \in
\dd} R \cdot \aaa$) for some $i$, then $\sg(\bfh) = 0$.

\item If this does not happen, i.e., $\Pi = (\Z / 2\Z)^k$ for some $k$,
and the hypotheses of the previous part hold, then 
\[ \sg(\aaa_s^n) = \sum_{\substack{l_i > 0\ \forall i \\ l_1 + \dots +
l_k = n}} \binom{n}{l_1, \dots, l_k} \prod_{i=1}^k \g_i(\aaa_s)^{l_i}, \]

\noindent where the $\g_i$'s are any set of generators for $\Pi$. In
particular, this vanishes if $k>r$, where $\sum_{j=1}^r 2^{s_j}$ is the
binary expansion of $n$.
\end{enumerate}
\end{theorem}

\begin{remark}\hfill
\begin{enumerate}
\item {\bf Warning.} One should not confuse the $h_i$'s here with
elements of $\mfh$; indeed, $h_i \in A$, so they really are in $\mfh^*$.

\item The coefficient above is just the multinomial coefficient $n! /
(\prod_i l_i!)$, which we also denote by $\nomial{n}$, just as
$\binom{n}{k,n-k} = \binom{n}{k}$.
The last line in the theorem follows because this coefficient is odd if
and only if ($r, s_j$ as above) we can partition $\{ 2^{s_j} : j \}$ into
$k$ nonempty subsets, and the $l_i$'s are precisely the sums of the
elements in the subsets. (This fails, for instance, if some two $l_i$'s
are equal, or $k>r$.)
In turn, this fact follows (inductively) from the following easy-to-prove

\begin{lemma}\label{Lnom}
Suppose $p>0$ is prime, $p^s \leq n < p^{s+1}$ for some $s \geq 0$, and
$l_k \geq l_i\ \forall i$. If $l_k < p^s$ then $p$ divides
$\nomial{n}$. Otherwise $p$ divides neither or both of $\nomial{n}$ and
$\nomial{n-p^s}$.
\end{lemma}
\comment{
\begin{proof}[Proof of Lemma \ref{Lnom}]
For the first part, suppose $l_k = \max_i l_i$ is less than $p^s$. We
show that the multinomial coefficient is divisible by $p$. First note
that the multinomial coefficient is a product of binomial coefficients:
\[ \binom{n}{l_1, \dots, l_k} = \binom{l_1}{l_1} \binom{l_1 + l_2}{l_2}
\dots \binom{l_1 + \dots + l_i}{l_i} \dots \binom{n}{l_k} \]

\noindent and each of the factors is an integer. Let us suppose that $l_1
+ \dots + l_{i-1} < p^s \leq l_1 + \dots + l_i$. One can then show that
the corresponding binomial coefficient factor on the right-hand side
above, is divisible by $p$, if $l_r < p^s\ \forall r$.\\

The second part follows if one shows that $p \nmid \displaystyle \frac{n
\dots (n-p^s+1)}{l_k \dots (l_k - p^s+1)}$, i.e., that the factors of $p$
in both the numerator and denominator are the same in number. But this is
clear, because by the assumption, there is no multiple of $p^{s+1}$ that
occurs in either expansion, so modulo $p^{s+1}$, both sets of numbers are
equivalent to $1, 2, \dots, p^s$. Therefore both products (or the powers
of $p$ that exactly divide them) equal the power of $p$ that divides
$(p^s)!$, and this common number is $p^{s-1} + p^{s-2} + \dots + 1 =
\frac{p^s - 1}{p-1}$.

Why does the second part follow now? Well, we have the easy-to-verify
equation:
\[ \nomial{n} = \frac{n \dots (n-p^s+1)}{l_k \dots (l_k - p^s + 1)}
\nomial{n-p^s} \]
\end{proof}
}
\end{enumerate}
\end{remark}

The rest of the subsection is devoted to the proof of the theorem. We
once again mention a result crucial to the proof, then use it to prove
the theorem, and conclude by proving the key claim.\medskip

\noindent {\bf Key claim.}
($\ch(R)$ arbitrary.) If $W$ contains a Dynkin subgraph $\Omega$ of type
$A_2$ or $G_2$, then both the simple roots in $\Omega$ are killed by all
$\la \in \mfh^W$. If $\Omega$ is of type $B_2$, then the long root in
$\Omega$ is killed by all $\la$.\medskip

\begin{proof}[Proof modulo the key claim]
We now show the theorem.
\begin{enumerate}
\item First suppose that $\ch(R) \neq 2$. If $\la \in \mfh^W$, then
$\la(\aaa) = \la(s_{\aaa}(\aaa)) = -\la(\aaa)$, whence $\la(\aaa) = 0$
for all $\aaa \in \dd$, and $\mfh^W = 0$.\medskip

For the other claims, use the classification of simple Lie algebras in
terms of Dynkin diagrams, as mentioned in \cite[Chapter 3]{Hum1}. To show
that a weight $\la$ kills all of $\mfh^*$, it suffices to show that
$\la(\aaa) = 0\ \forall \aaa \in \dd$, i.e., that it kills each simple
root, or node of the corresponding Dynkin diagram.

If the Dynkin diagram of a Lie algebra has (a sub-diagram of) type $A_2$
or $G_2$, then both nodes of that diagram (or both $\aaa_i$'s) are killed
by all weights $\la \in \Gamma$, by the key claim above. This automatically
eliminates all diagrams of type $A_n$ for $n>1$, as well as all
$D,E,F,G$-type diagrams, leaving only type $A_1$ among these.

Moreover, for types $B,C$, at most one simple root (the ``last" one) is
not killed by all $\la$'s. If this root is long, then it is also killed
by the key claim above (as a part of a $B_2$), and we are done.

\item First, $\la(s_\aaa^2) = \la(s_\aaa)^2 = 1$, whence $\la(s_\aaa) =
\pm 1 = 1\ \forall \aaa \in \dd$, if $\ch(R) = 2$. This implies that
$\la(w) = 1 = \vi(w)$ for all $w \in W, \la \in \Gamma$.
Next, Theorem \ref{T6} above tells us that if $\Pi$ has an element of
order 4, then $\sg(\bfh) = 0$.
Finally, if some $h_i$ has no ``$\aaa_s$-contribution", then it is killed
by all $\la$, by the previous part, so $\la(\bfh) = 0\ \forall \la \in
\Gamma$.

\item As we remarked after Theorem \ref{T8}, $\Pi = (\Z / 2\Z)^k$ in
characteristic 2, if $\Pi$ does not have an element of order 4. (Reason:
$\Gamma = \{ \vi \} \times \mfh^W \cong (\mfh^W, +)$ is a free $R$-module
by the previous part, and $2 \Gamma = 0$.)

We now perform the computation. For this, suppose that $h_i - c_i \aaa_s$
is in the $R$-span of $\{ \aaa \in \dd : \aaa \neq \aaa_s \}$ (note that
in the case of $A_1$, the condition $h_i \in R \cdot \aaa_s$ is
automatic). Then $\sg(\bfh) = \left( \prod_i c_i \right) \cdot
\sg(\aaa_s^n)$, so it suffices to compute $\sg(\aaa_s^n)$.

If $\{ \g_i \}$ is any set of generators (or $\Z/2\Z$-basis) for $\Pi$,
then
\comment{
we now compute:
\[ \sg(\aaa_s^n) = \sum_{I \subset [k]} \g_I(\aaa_s)^n \]

\noindent As in the key claim that was used to prove Theorem \ref{T8}
above, all monomials are of the form $\prod_{j \in I}
\g_j(\aaa_s^{l_j})$, where $I \subset [k]$ and $l_j > 0, \sum_j l_j =
n$. Moreover, the coefficient of such a monomial is a multiple of 2 (and
hence vanishes) if $I \neq [k]$.

It thus remains to compute the coefficient for the cases when $I = [k]$,
and this is obtained precisely as a multinomial coefficient by
definition, since it occurs exactly so many times in only one term,
namely, $\g_{[k]}(\aaa_s)^n$.
}
the desired equation actually holds if we sum over all {\it nonnegative}
tuples $l_i$, that add up to $k$. Thus, the proof is similar to that of
the key claim used to prove Theorem \ref{T8} above; simply note that if
$I \subsetneq [k]$, then every $\prod_{j \in I} \g_j(\aaa_s^{|I_j|})$
occurs with an even coefficient.
\end{enumerate}
\end{proof}

\noindent Finally, we prove the key claim.

\begin{proof}[Proof of the key claim]
It helps to look at the pictures of these rank 2 root systems (drawn in
\cite[Chapter 3]{Hum1}). We use the $W$-invariance of $\la|_V\ \forall
\la \in \Gamma$.

Consider the system $A_2$, with simple roots $\aaa, \beta$. Given $\la
\in \Gamma$, $\la(\aaa) = \la(\beta) = \la(\aaa + \beta)$, whence
$\la(\aaa) = \la(\beta) = 0$.

The root system $G_2$ has two subsystems of type $A_2$, whence each $\la$
must kill both subsystems.

Now consider $B_2$, with long root $\aaa$ and short root $\beta$.
Clearly, $\beta + \aaa$ is another short root, whence $\la(\beta + \aaa)
= \la(\beta)$, and we are done.
\end{proof}

\subsection{The reductive case}

We conclude by mentioning what happens in the reductive case. This uses
the results proved in the simple case above. Recall also that the
notation for this situation was set when we defined degenerate affine
Hecke algebras with trivial parameter earlier. This notation will be used
freely here, without recalling it from above.

Let $V'$ be the direct sum of $V_0$ and the $R$-span of all the unique
short simple roots $\aaa_{i,\ short}$ inside any of the simple components
$V_i = \mfh^*_i$ of the ``correct" type (not $G_2$). Let the other simple
roots in $\dd$ span the $R$-submodule $V''$. Then $V = V' \oplus V''$,
and each $\la \in \Gamma$ kills $V''$. There now are two cases.\medskip

\noindent {\bf Case 1.} $\ch(R) \neq 2$.
Then $\la$ in fact kills all $\aaa \in \dd$, because $\la(\aaa) =
\la(s_\aaa(\aaa)) = -\la(\aaa)$. This means that we are left with $V_0$,
i.e., if for all $i$, $h_i - v_{0,i} \in \oplus_{j>0} V_j$ for some
$v_{0,i} \in V_0$, then $\sg(\bfh) = \sg(\prod_i v_{0,i})$.

Next, recall that $\sg = \Gamma_W \times (V^*)^W$, so we are reduced to
the case of every $\la$ being represented (on $V_0$) by some element of
$V_0^* = (V_0^*)^W$. We conclude this case by noting that some (partial)
results on how to compute this were included in the previous
section.\medskip

\noindent {\bf Case 2.} $\ch(R) = 2$.
Then $\la(w) = 1$ for all $w, \la$, as seen above. Moreover, we are left
only to consider the case of all $h_i \in V'$. Now, $\Gamma = \vi_W
\times (V')^*$, whence any finite subgroup $\Pi = (\Z / 2\Z)^k$ for some
$k$ (since it too is a $\Z / 2\Z$-vector space). In this situation,
Theorems \ref{T6} and \ref{T8} (and \ref{T7} as well) give us some
information on how to compute $\sg(\bfh)$.

\section{An example that attains any value}

We conclude with examples where $\sg(\bfh)$ can take any value in $R$, if
the $h_i$'s are merely skew-primitive.\medskip

\noindent {\bf Example 12} (A skew-primitively generated algebra){\bf .}
By Proposition \ref{P7} above, if all $h_i$'s are pseudo-primitive, then
$\sg(\bfh) = 0$ if $\ch(R) \nmid |\Pi|$ - whereas if $\ch(R)$ divides
$|\Pi|$, then this case was analyzed in Section \ref{Spseudo} above.

One can ask if such results hold in general, i.e., for products of
skew-primitive elements. (Note by Theorem \ref{Tskew} that we need:
$\ch(R) \nmid |\Pi|$.)

For the example that we now mention (for groups $\Pi$ of {\em even}
order), one needs to {\bf assume} the following:
\begin{enumerate}
\item $\ch(R) > 2^n$ and $\exp(\Pi)$, or $\ch(R) = 0$ and $R \supset \Q$;
and
\item If $d = \exp(\Pi)$ is the exponent, then $d$ is even, and there
exists a primitive $d$th root of unity in $R$, say $z$.
\end{enumerate}

\noindent Beyond this, given $n,\Pi$ (of even order), and $r \in R$, we
will produce the desired Hopf algebra $\scrh$, a group of weights $\Pi
\subset \Gamma_{\scrh}$, and skew-primitive $h_1, \dots, h_n \in H$, such
that $\sg(\bfh) = r \in R$.\medskip

Given $\Pi$, use the Structure Theorem for finite Abelian groups, to
write: $\Pi = \bigoplus_{i=1}^k (\Z / d_i \Z)$, with $d_1 | d_2 | \dots |
d_k = \exp(\Pi)$. Then $d_k$ is even, since $\Pi$ has even order. Now
define $\scrh$ to be the commutative $R$-algebra freely generated as:
$\scrh = R[R^n] \otimes R[\Z^k]$. In other words, $R$ is generated by
$h_1, \dots, h_n, g_1^{\pm 1}, \dots, g_k^{\pm 1}$, with the relation
that they all commute (and that the $g_i$'s are invertible).

Now define the $g_i$'s to be grouplike and $\dd(h_j) = g_k \otimes h_j +
h_j \otimes 1$. Also define (for all $i,j$):
\[ \vi(g_i) = 1,\ S(g_i) = g_i^{-1},\ \vi(h_j) = 0,\ S(h_j) = -g_k^{-1}
h_j. \]

\noindent Since $\scrh$ is freely generated, the set of weights of
$\scrh$ is $R^n \times (R^\times)^k$. Since (it can be checked that)
$\scrh$ is also a Hopf algebra, the group operation is:
\begin{eqnarray*}
&& (a_1, \dots, a_n, z_1, \dots, z_k) * (a'_1, \dots, a'_n, z'_1, \dots,
z'_k)\\
& = & (a_1 + z_k a'_1, a_2 + z_k a'_2, \dots, a_n + z_k a'_n, z_1 z'_1,
\dots, z_n, z'_n).
\end{eqnarray*}

We now produce the desired example. Define $\g_j \in \Gamma_{\scrh}$ on
generators by: $\g_j(g_i) = z^{\delta_{ij} d_k / d_j}$, and $\g_j(h_i) =
0$ unless $j=k$. Moreover, $\g_k(h_i) = 1$ for $i<k$, and $\g_k(h_k) = (1
- z)^n |\Pi|^{-1} r'$ for some $r' \in R$ (which we define later, and
which depends on $n$).

It is now easy to check that each $\g_j$ is of order $d_j$, and the
$\g_j$'s generate a subgroup of $\Gamma_{\scrh}$ isomorphic to $\Pi$.
Moreover, $\g_1, \dots, \g_{k-1}$ all kill $h_1, h_2, \dots, h_n$.
Adopting the notation of Proposition \ref{Pstab}, $\g_j \in \Gamma_{h_i}$
for $j < k$ and all $i$. Since $\Gamma_{h_i}$ is a subgroup of
$\Gamma_{\scrh}$, hence $\Pi_1 := \langle \g_1, \dots, \g_{k-1} \rangle
\subset \Gamma_{h_i}$ for all $i$; by Proposition \ref{Pstab}, $\Pi_1
\subset \Gamma_{\bfh}$.\medskip

Now use Lemma \ref{L1} (noting that $\Pi_1$ is normal in the abelian
group $\Pi$); then $\sg(\bfh) = |\Pi / \Pi_1| \Sigma_{\Pi_1}(\bfh)$. Use
Theorem \ref{Tabel}: $f_i = \frac{\g_k(h_i)}{\g_k(g_k) - 1}$, so
\[ (-1)^n \prod_i f_i = (-1)^n \frac{(1 - z)^n |\Pi|^{-1}
r'}{\prod_{i=1}^n (z - 1)} = \frac{r'}{|\Pi|}. \]

Moreover, $S = \{ I \subset \{ 1, 2, \dots, n \} : g_k^{|I|}$ is fixed by
$\Pi_1 \}$, i.e., all subsets $I$ such that $d_k | |I|$. Since $d_k$ is
even, this means that $(-1)^{|I|} = 1\ \forall I \in S$, whence
\[ \sg(\bfh) = |\Pi / \Pi_1| \Sigma_{\Pi_1}(\bfh) = \frac{|\Pi|}{|\Pi_1|}
\cdot |\Pi_1| \cdot \frac{r'}{|\Pi|} \cdot \sum_{m \geq 0} \binom{n}{m
d_k} = r' \sum_{m \geq 0} \binom{n}{m d_k}. \]

\noindent By assumption (on $R$), the summation is a unit in $R$, so
choosing $r'$ suitably, one obtains any $r \in R$ as our
answer.\qed\medskip

Also note (e.g., by \cite[Exercise 38, \S 1.2.6]{Knu}), that the
summation equals $\frac{1}{d_k} \sum_{l=0}^{d_k-1} (1 + z^l)^n$.\bigskip

\section*{Acknowledgements}

I thank Susan Montgomery and Nicol\'as Andruskiewitsch for their comments
and suggestions after reading a preliminary draft of this manuscript. I
also thank the referees for their feedback, which helped improve the
exposition of this paper.




\end{document}